%% file: Survey.tex
 \documentclass[final,5p,times,twocolumn,authoryear]{elsarticle}

\usepackage{amssymb}

\journal{Annual Reviews in Control}

\usepackage{graphicx}
\usepackage{natbib}
\usepackage{mathtools}
\usepackage{tikz}
\usepackage{amsmath,amssymb,amsfonts}

\newtheorem{assu}{Assumption}

\newtheorem{rmk}{Remark}

\def\TM{\textcolor{black}}

\begin{document}

\begin{frontmatter}



\title{Guarantees for data-driven control of nonlinear systems using semidefinite programming: A survey \tnoteref{t1}}


\author[Stuttgart]{Tim Martin\corref{cor1}} \ead{tim.martin@ist.uni-stuttgart.de}
\author[Uppsala]{Thomas B. Schön} \ead{thomas.schon@it.uu.se}
\author[Stuttgart]{Frank Allgöwer} \ead{frank.allgower@ist.uni-stuttgart.de}

\affiliation[Stuttgart]{organization={University of Stuttgart, Institute for Systems Theory and Automatic Control}, country={Germany}}
\affiliation[Uppsala]{organization={Uppsala University, Department of Information Technology}, country={Sweden}}

\tnotetext[t1]{Frank Allgöwer thanks the funding by the Deutsche Forschungsgemeinschaft (DFG, German Research Foundation) under Germany's Excellence Strategy - EXC 2075 - 390740016 and under grant 468094890. The work of Thomas B. Schön is funded by \emph{Kjell och Märta Beijer Foundation} and by the project \emph{NewLEADS - New Directions in Learning Dynamical Systems} (contract number: 621-2016-06079), funded by the Swedish Research Council. Tim Martin and Frank Allgöwer also acknowledge the support by the Stuttgart Center for Simulation Science (SimTech).\\
	\copyright2023. This manuscript has been accepted to Annual Reviews in Control ({https://doi.org/10.1016/j.arcontrol.2023.100911}) and is made available under a Creative Commons Licence CC-BY-NC-ND.}
\cortext[cor1]{Corresponding author}  

\begin{abstract}
This survey presents recent research on determining control-theoretic properties and designing controllers with rigorous guarantees using semidefinite programming and for nonlinear systems for which no mathematical models but measured trajectories are available. Data-driven control techniques have been developed to circumvent a time-consuming modelling by first principles and because of the increasing availability of data. Recently, this research field has gained increased attention by the application of Willems' fundamental lemma, which provides a fertile ground for the development of data-driven control schemes with guarantees for linear time-invariant systems. While the fundamental lemma can be generalized to further system classes, there does not exist a comparable data-based system representation for nonlinear systems. At the same time, nonlinear systems constitute the majority of practical systems. Moreover, they include additional challenges such as \TM{data-based surrogate models that prevent system analysis and controller design by convex optimization.} Therefore, a variety of data-driven control approaches has been developed with different required prior insights into the system to ensure a guaranteed inference. In this survey, we will discuss developments in the context of data-driven control for nonlinear systems. In particular, we will focus on methods based on system representations providing guarantees from finite data, while the analysis and the controller design boil down to convex optimization problems given as semidefinite programming. Thus, these approaches achieve reasonable advances compared to the state-of-the-art system analysis and controller design by models from system identification. Specifically, the paper covers system representations based on extensions of Willems' fundamental lemma, set membership, kernel \TM{techniques}, the Koopman operator, and feedback linearization.
\end{abstract}



\begin{keyword}
Data-driven control \sep Data-driven system analysis \sep Nonlinear systems \sep Semidefinite programming 


\end{keyword}

\end{frontmatter}



\input{01_Intro}

\input{02_Preview}

\input{03_Polynomial_Approx}

\input{04_GP}
\input{05_LPV}

\input{06_Koopman}

\input{07_Feedback_linearization}

\input{08_Conclusion}




\end{document}

%% file: 01_Intro.tex
\section{Introduction}\label{Sec_Intro}

Model-based control techniques suppose the access to a mathematical model that describes the behavior of a system over time. The description of the dynamics can be given as difference equations in discrete time or as differential equations in continuous time of the system's inputs and states or outputs. 
Besides the synthesis of a controller to influence the behavior of the system, system analysis aims to provide valuable insights into the system by the verification of control-theoretic properties such as dissipativity \citep{DissiWillems}. Subsequently, these properties can be used for a controller design using feedback laws \citep{Khalil}, e.g., the small-gain theorem or the interconnection of passive systems.

One possible derivation of models is based on first principles, for instance, Newton's laws of motion, Kirchhoff's circuit laws, and the first and second law of thermodynamics. However, their application often requires expert knowledge, calls for a priori simplifications to obtain suitable models for control, or is more time consuming than the controller design itself. At the same time, the goal intrinsically is the controller instead of a model. 

For these reasons, data-driven approaches \citep{DDC_Survey} have gained in popularity. There system properties are verified and controllers designed from measured trajectories of the underlying system. System identification \citep{SysID,SysId} represents a so-called indirect data-driven method because first a model is identified from data and then analyzed or a controller is derived by model-based techniques. However, here the mismatch between the identified model and the underlying system is often obscure. Indeed, the investigation of the model mismatch is even for the identification of linear time-invariant (LTI) systems an active research field \citep{Oymak}. \TM{At the same time, if this mismatch is unknown, then} the design of a controller with closed-loop stability and performance guarantees is jeopardized though inherent guarantees of the control design procedure.  

Investigations of direct data-driven techniques for LTI systems without an intermediate modelling step include PID control \citep{PID}, adaptive control \citep{LinAdaptive}, iterative feedback tuning \citep{IDT}, virtual reference feedback tuning \citep{VRFT}, reinforcement learning \citep{Reinforcement}, unfalsified control \citep{UnfalContr}, and subspace-based LQG-control \citep{SubSpaceID}. 

Moreover, the renewed interest in the behavioral approach and the fundamental lemma from \cite{Willems} in the context of data-driven control has led to a framework \citep{Markovsky_Survey} for data-driven system analysis and various control schemes. Furthermore, inspired by Willems' fundamental lemma, \cite{PersisLinear} presents a parametrization of the closed loop of an LTI system and a state feedback based on data matrices. This result also has led to further data-driven control synthesis approaches \citep{Markovsky_Survey}.  

While the data for Willems' fundamental lemma would allow the exact identification of the LTI system, the data-informativity framework of \cite{Informativity,Inform_Survey} examines the question, when data are informative enough to draw conclusions concerning controllability, stabilizability, etc. Therefore, the amount of data that is needed to identify a system is in general larger than what is needed to control the system. The framework is based on the set of all LTI systems explaining the data. This resembles a set-membership approach \citep{Set_memerbship_Lin}, which is, e.g., exploited in \cite{Superstable} to design a controller for superstability from noisy input-output data by linear programming. 

These recent developments have established for LTI systems a comprising framework for data-driven system analysis and control including, e.g., verification of various system properties, optimal control, robust control, and predictive control. However, analogous results are missing for nonlinear systems due to the manifold of additional challenges. Indeed, the identification of a nonlinear system only from a finite set of samples is not possible. Instead, additional a priori insights into the system are required. 
Moreover, although many controller design techniques have inherent closed-loop stability guarantees, they are jeopardized due to the not exactly known nonlinearities. Thus, no end-to-end guarantees for the closed loop can be recovered. Additionally, estimating a bound between the nonlinearity and its \TM{estimation} from finite data is \TM{nontrivial}. Lastly, a direct analysis of nonlinear systems leads in general to nonconvex optimization problems, for example, the estimation of the region of attraction.   

Data-driven approaches for nonlinear systems with guarantees include nonlinear adaptive control \citep{Adaptive}, learning-based model predictive control (MPC) \citep{Learn_MPC}, reinforcement learning with safety guarantees \citep{Reinf_NL}, and stability verification \citep{Scenario} using scenario optimization \citep{Scenario_opt}. Furthermore, neural networks are exemplarily applied in \cite{NN_controller} to simultaneously learn the dynamics, controller, and a Lyapunov function to provide stability guarantees under a known approximation error of the neural network. Under a known upper bound on the Lipschitz constant of the system's dynamics, set membership \citep{Milanese} and Kinky inference \citep{Kinky} provide a framework for nonlinear systems to design a controller, for instance, by online prediction \citep{Set_membership_control} or online certificate function control \citep{Lipschitz_online_control}. Related to this framework, \cite{Montenbruck} presents a data-driven system analysis via the input-output mapping of a nonlinear system and \cite{MartinIterative} examines an extension by an iterative sampling scheme.   

This article focuses on data-driven methods based on system representations providing rigorous guarantees and enabling verification of control-theoretic properties and the design of controllers via semidefinite programming (SDP). To circumvent \TM{nonconvex} optimization though nonlinear system dynamics, convex relaxations and various linearization and polynomialization methods are employed. Due to noisy measurements or the approximation of the system dynamics, the underlying dynamics can not exactly be identified. Nevertheless, the presented approaches provide guarantees by first obtaining a set of systems consistent with the data and a bound on the approximation error. Combined with robust control techniques and a convex optimization via SDPs, rigorous end-to-end guarantees for the determined system properties and closed loop are ensured from finite data. The motivation for the development of these methods is explained in the following points (i)-(vi).

(i) The presented methods \TM{derive} system representations suitable for data-driven system analysis and the design of various controller schemes from data. \TM{Therefore, these methods strive to} establish a framework for nonlinear data-driven system analysis and control, analogously to the frameworks for LTI systems by the fundamental lemma \citep{Willems} and data informativity \citep{Informativity}. 

\TM{(ii)} System identification, e.g., based on neural networks or Gaussian processes, can result in precise but strong nonlinear, models. Whereas a system analysis or a state feedback design based on these models would be a nonconvex optimization problem, the presented system representations are tailored for solving many control-related problems by SDPs. 
 
\TM{(iii)} As shown in Section~\ref{Sec_SDP}, SDPs can be solved in practice often in a tractable way by relying on well-established algorithms and solvers. 

\TM{(iv)} While system identification aims to approximate the dynamics as precise as possible, the system representations presented here are motivated to verify system properties or design a controller. As in the LTI case \citep{Informativity}, one expect that the identification requires more data than the verification or controller synthesis problem.

\TM{(v)} System identification techniques for nonlinear systems not always provide error bounds to ensure guarantees.

(vi) The presented system representations leverage, among others, a set of systems feasible with the observed data, similar to a set-membership approach \citep{Milanese}. While the existing literature considers Lipschitz approximations for specific control structures and data-inefficient system analysis \citep{Montenbruck}, the presented methods exploits more general robust control techniques to, e.g., include performance criteria for more general closed-loop structures. In contrast to set-membership identification \citep{Milanese_SysId}, where the feasible system set is leveraged to obtain a model and its mismatch, the methods here exploits the system set directly for control.

The survey is organized as follows. We begin with the introduction of some notation and a motivation of SDPs in control in Section~\ref{Sec_SDP}. Afterwards, we will briefly report the presented data-driven approaches in Section~\ref{Preview} and then provide a more detailed discussion of their key ideas. Therefore, we will focus on the data-based representation of the nonlinear system rather than their application to specific control problems. Section~\ref{Sec_Discussion} shows further discussion of the presented methods and Section~\ref{Sec_Conclusion} concludes the article. 

\section{Notation}\label{Sec_Notation}

Throughout the article, we denote the set of natural numbers by $\mathbb{N}$, the natural numbers including zero by $\mathbb{N}_0$, and the set of real numbers by $\mathbb{R}$. The Euclidean norm of a vector $v\in\mathbb{R}^n$ is denoted by $||v||_2$. Furthermore, $I$ and $0$ corresponds to the identity and the zero matrix of suitable dimensions, respectively. The Frobenius norm of a matrix $M\in\mathbb{R}^{n\times m}$ is denoted by $||M||_\text{Fr}$. The right inverse of a full-row-rank matrix $R\in\mathbb{R}^{n\times m}$ is denoted by $R^\dagger\in\mathbb{R}^{m\times n}$. For a symmetric matrix $A=A^T$, $A\succ0$ or $A\succeq0$ denote that $A$ is positive definite or positive semidefinite, respectively. Analogously, $A\prec0$ or $A\preceq0$ if $A$ is negative definite or negative semidefinite, respectively.

\section{Semidefinite programming in control}\label{Sec_SDP}

\begin{table*}
	\begin{center}
		\begin{tabular}{c|c|c}\label{OverviewTable}
			Section & Specific methods & Main references\\\hline\hline
			Polynomial approximation & Polynomial interpolation &\cite{MartinTP1,MartinTP2}, \cite{DePersis2},\\&&\cite{MartinCS} \\
			& Polynomial subclass& see\footnotemark[1] \\
			& Data-based closed-loop description& \cite{PersisLinear}, \cite{DePersis1}, \\&&\cite{DDC_orbit}, \cite{Persis_quad_constr}\\\hline
			Kernel regression &Linear sector from GP & \cite{Fiedler}\\
			& Linearized kernel & \cite{GP_RobCon}, \cite{Persis_kernel}\\
			&Polynomial kernel & \cite{Poly_Kernel}\\\hline
			LPV embedding &LPV system &\cite{LPV_FundamentalLemma}, \cite{LPV_direct},\\&& \cite{LPV_set_membership}, \cite{LPV_SetMem2}\\
			& LPV embedding & \cite{LPV_nonlinear}, \cite{LPV_nonlinear2}\\
			& Extended linearization & \cite{Sznaier} \\\hline
			State lifting & Koopman & \cite{Koopman_Guarantees}, \cite{Strasser_Koopman}\\\hline
			Feedback linearization & Nonlinearity cancellation & \cite{NonlinearityCancell}	\\
			& Flat system & \cite{FeedbackLin2,FeedbackLin,FeedbackLin3,ExpDesign_FeebackLin,MPC_FeebackLin}	 
		\end{tabular}
	\caption{Organization of presented methods and references.}
	\end{center}	
\end{table*}

This section motivates the application of SDPs in system analysis and control. To this end, we first provide a brief introduction of SDPs following \cite{SDP_intro}. Second, we comment on the usefulness of SDPs in control.

While SDPs are introduced in different forms, we remain on a control perspective. There an SDP minimizes a linear function subject to a constraint given by the definiteness of an affine function of symmetric matrices. More precisely,
\begin{equation}\label{SDP}
\begin{aligned}
	\min_{x\in\mathbb{R}^{n}}\hspace{0.2cm} &c^Tx\\
	\text{subject to } &F(x)=F_0+\sum_{i=1}^{n}x_iF_i\succeq0,
\end{aligned}
\end{equation}
with $c\in\mathbb{R}^{n}$ and symmetric matrices $F_i\in\mathbb{R}^{m\times m},i=0,\dots,n$. Since the objective function as well as $F(x)$ are linear in the optimization variable $x$, optimization problem \eqref{SDP} is convex. In this case, all local optima are globally optimal such that convex optimization problems are theoretically tractable.

\begin{figure}
\centering
\begin{tikzpicture}
\node[inner sep=0pt] (p1) at (0,0) {\includegraphics[width=.6\textwidth]{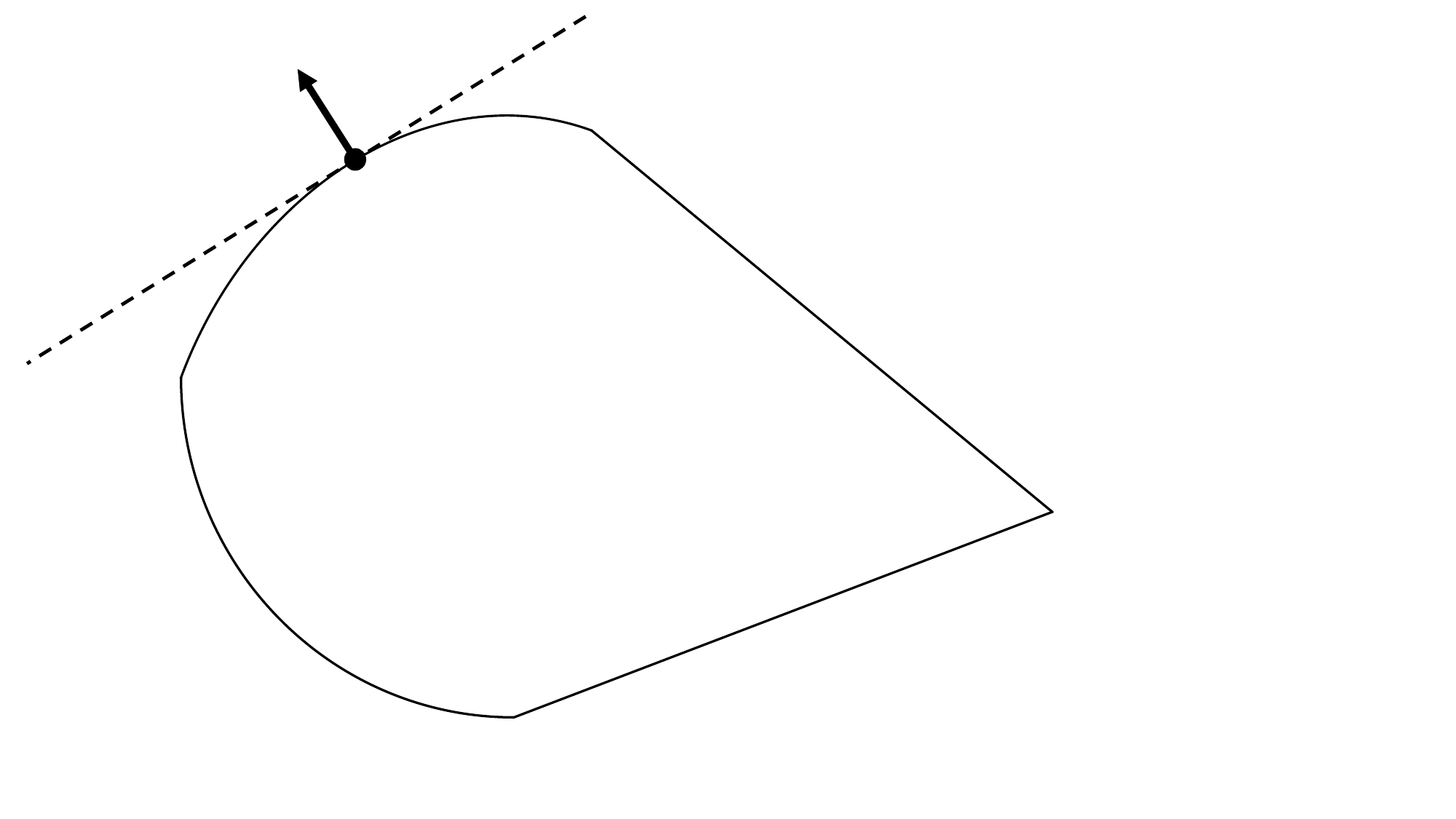}};

\node[inner sep=0pt] (p1) at (-1,0) {\large $F(x)\succeq0$};
\node[inner sep=0pt] (p1) at (1.5,2) {\large $F(x)\nsucceq0$};
\node[inner sep=0pt] (p1) at (-3,2.7) {\large$c$};
\node[inner sep=0pt] (p1) at (-2.5,1.6) {\large$x_\text{opt}$};

\end{tikzpicture}
\caption{Graphical illustration of an SDP.}
\label{Fig.SDP}
\end{figure}
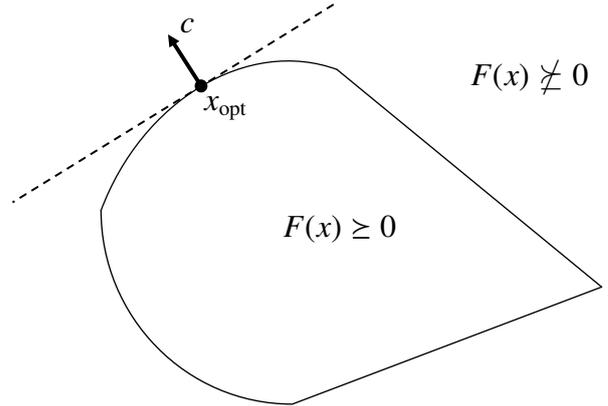

The optimization problem of SDP \eqref{SDP} can be illustrated as in Figure~\ref{Fig.SDP}. To find its optimal solution, we need to push the dotted line as far as possible into the direction of $-c$ while not intersecting the feasible set $F(x)\succeq0$. As illustrated, the feasible set is convex and the optimal solution $x_\text{opt}$ lies in general on its boundary, i.e., the matrix $F(x_\text{opt})$ is singular.

Besides the theoretical tractability due to convexity, SDPs are attractive as many problems from combinatorial optimization and control theory can be recast as an SDP. Indeed, SDPs can handle constraints given by linear matrix inequalities (LMI) $F(x)\succeq0$, convex quadratic inequalities, lower bounds on matrix norms, lower bounds on determinants of positive semidefinite matrices, and polynomial inequalities via sum-of-squares (SOS) hierarchies (see Section~\ref{Sec_SOS}).

Furthermore, SDPs generalize linear programming with matrix inequalities instead of componentwise inequalities between vectors. By this connection, many results and algorithms from linear programming extend to SDPs even though the latter is more general. For instance, most interior-point methods for linear programming can be generalized to SDPs. To this end, barrier functions are introduced which tend to infinity as points approach the boundary of the feasible set. Thereby, the constraint optimization \eqref{SDP} can be reformulated into an unconstrained one, which can be solved efficiently by Newton iteration techniques. Similar to linear programming, these iterations have polynomial worst-case complexity and perform very well in practice. For instance, \cite{SDP_intro} provides the rule-of-thumb that interior-point methods solve SDPs in $5-50$ iterations, where each iteration corresponds to a least-squares problem of the same size as the original problem.

We conclude this section by a simple application of SDPs in a control context. For that purpose, we want to check whether an LTI system $\dot{x}(t)=Ax(t)+Bu(t),y(t)=Cx(t)+Du(t)$ with state $x(t)\in\mathbb{R}^{n_x}$, input $u(t)\in\mathbb{R}^{m}$, and output $y(t)\in\mathbb{R}^m$, is passive. According to \cite{DissiWillems}, we need to search for a positive definite function $S(x)$ such that
$\dot{S}(x(t))\leq y(t)^Tu(t)$ for all $x(t)\in\mathbb{R}^{n_x}$ and $u(t)\in\mathbb{R}^{m}$. For the quadratic ansatz $S(x)=x^TPx,P\succeq0$, we obtain the condition
\begin{equation}\label{Pass_cond}
	\begin{bmatrix}x\\u\end{bmatrix}^T\underbrace{\begin{bmatrix}-PA-A^TP & -PB+\frac{1}{2}C\\ -B^TP+\frac{1}{2}C & D\end{bmatrix}}_{=L(P)}\begin{bmatrix}x\\u\end{bmatrix}\geq0,
\end{equation}
which is implied by $L(P)\succeq0$. By introducing each element of $P=\begin{bmatrix}p_{11} & \cdots\\ \vdots & \ddots\end{bmatrix}$, one can see that $L(P)\succeq0$ is an LMI with the optimization variables $p_{ij},i=1,\dots,n_x,j=1,\dots,n_x$. Hence, we can check for passivity of an LTI system by the SDP
\begin{align*}
	\min_{P\in\mathbb{R}^{n_x\times n_x}}\hspace{0.2cm} &0\\
	\text{subject to } &\begin{bmatrix}P & 0\\ 0 & L(P)\end{bmatrix}\succeq0.
\end{align*}

When replacing the LTI dynamics by a nonlinear $\dot{x}(t)=f(x(t),u(t)),y(t)=h(x(t),u(t))$, then condition \eqref{Pass_cond} results into a nonlinear matrix inequality. The resulting optimization problem is therefore not an SDP in general. \TM{To rely on the well-established solvers for SDPs,} the methods presented here not only need to infer on the dynamics of the system but also derive a system description suitable for a system analysis and a controller design by SDPs.


\subsection{Sum-of-squares optimization}\label{Sec_SOS}

Many control-related problems, e.g., the verification of Lyapunov stability of a polynomial system, include polynomial inequality constraints. To solve these NP-hard problems, polynomial inequalities can be relaxed by SOS decomposition leading to an SDP with LMI constraints. We will briefly introduce SOS optimization here as various presented methods rely on this relaxation. 

Consider a real polynomial in $x=\begin{bmatrix}x_1 & \cdots &x_n\end{bmatrix}^T\in\mathbb{R}^n$ of degree $d$  
\begin{equation*}
p(x)=\sum_{\alpha\in\mathbb{N}_0^n,|\alpha|\leq d} a_\alpha x^\alpha,
\end{equation*}
with vectorial indices $\alpha^T=\begin{bmatrix}\alpha_1 & \cdots & \alpha_n\end{bmatrix}^T\in\mathbb{N}_0^n$, $|\alpha|=\alpha_1+\cdots+\alpha_n$, real coefficients $a_\alpha\in\mathbb{R}$, and monomials $x^\alpha=x_1^{\alpha_1}\cdots x_n^{\alpha_n}$. Let $\mathbb{R}[x]$, $\mathbb{R}[x]^m$, and $\mathbb{R}[x]^{m\times n}$ denote the set of all real polynomials, all $m$-dimensional polynomial vectors, and all $m\times n$ polynomial matrices, respectively. Then a matrix $P\in\mathbb{R}[x]^{n\times n}$ is an SOS matrix if there exists a matrix $Q\in\mathbb{R}[x]^{{m}\times {n}}$ such that $P(x)=Q(x)^TQ(x)$. According to the square matricial representation \citep{SOS_decomp}, the search for $Q$ boils down to checking the feasibility of an LMI. 

By the SOS decomposition $P(x)=Q(x)^TQ(x)$, all SOS matrices are positive semidefinite for all $x\in\mathbb{R}^{n}$. However, not all polynomial positive semidefinite matrices are SOS. Thus, an SOS condition corresponds to a relaxation of a positive semidefiniteness condition of a polynomial matrix. Nevertheless, SOS optimization is widely used in control for the following two reasons. The 
Positivstellensatz from \cite{Putinar} allows to check for non-negativity of a polynomial on a compact semialgebraic set. Together with a matrix-version of this result \citep{SchererSOS}, a regional analysis of control systems is possible. 
Moreover, the moment-SOS hierarchy \citep{SOS_hierachy} shows that the relaxation converges with increasing degree of polynomials. 

One drawback of SOS optimization is scalability. Indeed, for a polynomial $p(x)$ with degree $2d$, the SOS relaxation of $p(x)\leq0$ leads to $p(x)=z(x)Qz(x)$ and $Q\succeq0$ with $Q$ of dimensions $\begin{pmatrix}n+d\\d\end{pmatrix}\times\begin{pmatrix}n+d\\d\end{pmatrix}\approx n^d\times n^d$. If the scalability of standard SOS relaxation prevents the application to large problems, chordal sparsity \citep{SOS_Sparse} for the obtained SOS problem might improve scalability issues. Furthermore, SOS relaxation could be replaced by B-spline relaxations, which shows less conservatism and less computational demand for, e.g., LPV controller design with polynomial parameter dependence \citep{B_spline_relax}.

%% file: 02_Preview.tex
\section{Preview}\label{Preview}

The purpose of this section is to provide a preview of the data-driven methods covered in this survey. They are organized by the derivation technique for their data-based system representation as given in Table~1. \footnotetext[1]{\cite{MartinPoly}, \cite{MartinIQC}, \cite{DePersis1}, \cite{Petersen}, \cite{FarkasLem}, \cite{Poly_Invariance}, \cite{SafeControl2}, \cite{Strasser}, \cite{Berberich}, \cite{Findeisen}, \cite{Input_constr}}


%

\subsection{Data-based polynomial approximation}\label{Sec_Poly_Approx}

In Section~\ref{Sec_PolyApprox}, we will take a closer look at characterizations of the unknown nonlinear dynamics by data-based polynomial approximation. Polynomial approximation has been widely used to deal with nonlinear systems in control theory \citep{Carleman} and in application, e.g., by Taylor linearization of the system dynamics. Thus, this data-based system representation is intuitive from a control perspective. Furthermore, a polynomial representation allows for the verification of system properties and for the design of controllers by SOS optimization. Moreover, a polynomial approximation does not require knowledge of a function basis containing the system dynamics. At the same time, the literature on polynomial interpolation \citep{Hermite} provides well-investigated approximation errors. These are essential to infer a tight set membership for the nonlinear system from data and to provide guarantees. 

For Taylor polynomials (TP), \cite{MartinTP2} verifies dissipativity properties, \cite{MartinTP1} determines incremental dissipativity, \cite{DePersis2} derives locally asymptotically stabilizing controllers, and \cite{MartinCS} obtains state-feedback laws to render an equilibrium globally asymptotically stable while satisfying closed-loop performance criteria. Since TPs commonly provide local approximations, \cite{MartinTP1} combines multiple TPs to refine the data-driven inference. 

If the error of the polynomial approximation vanishes, then the special case of polynomial systems is obtained. Data-driven control for unknown polynomial systems using SOS relaxation include the verification of dissipativity \citep{MartinPoly} and integral quadratic constraints (IQC) \citep{MartinIQC} and the controller synthesis \citep{DePersis1,Petersen,FarkasLem,Poly_Invariance,SafeControl2}. Related to the polynomial case, \cite{Strasser} considers the controller synthesis for rational systems. Related to the polynomial approximation by TPs, \cite{Berberich, Findeisen,Input_constr} investigate Lur'e-type systems \citep{Khalil} (Chapter 10.1) with a data-driven inference of the LTI part of the dynamics while assuming measurements and a known sector bound on the nonlinear part. Note that all these results excessively exploit techniques from robust control as linear fractional representation (LFR), Petersen's lemma, a matrix S-lemma \citep{vanWaarde}, or Farkas' lemma to provide guarantees though the system dynamics is not precisely known. 

In contrast to the previous set-membership approaches, \cite{PersisLinear} uses Taylor linearization to provide a data-driven representation of the linear part of the closed loop. Thereby, an equilibrium point can be rendered locally asymptotically stable by solving an optimization problem with LMI constraints. Moreover, the closed-loop characterization is extended to polynomial systems \citep{DePersis1}, periodic orbits \citep{DDC_orbit}, and Lur'e-type systems \citep{Persis_quad_constr}.

\subsection{Gaussian processes and kernel ridge regression}

Gaussian processes (GP) and kernel ridge regression constitute a flexible framework to approximate nonlinear functions in machine learning and nonlinear dynamics in system identification. Both regression methods provide the possibility to include prior knowledge and inherent uncertainty measures to derive guarantees for data-driven control. However, the obtained system representation is often strongly nonlinear due to nonlinear kernel functions. To deal with this nonlinearity, \cite{Umlauft} presents a controller design by feedback linearization and \cite{GP_Backstepping} by backstepping. Nonetheless, both require a certain structure of the system dynamics. Therefore, we will study in Section~\ref{Sec_GP} the following three approaches to achieve a controller synthesis by SDPs. 

\cite{Fiedler} learns a linear sector for the nonlinear parts of the dynamics from a GP to apply linear robust control afterwards. To this end, \cite{Fiedler2} establishes a statistical bound between the underlying nonlinear dynamics and the mean function of the GP. 

Instead of bounding the mismatch of the kernel regression by a sector, \cite{GP_RobCon} directly computes the Taylor linearization of the nonlinear GP-model around an equilibrium point for a linear robust controller design. Alternatively, \cite{Persis_kernel} suggests to stabilize the linear part of a kernel regression, while approximately cancelling its nonlinearity as in \cite{NonlinearityCancell}.  

\cite{Poly_Kernel} proposes to use polynomial kernels yielding a polynomial regression model and a polynomial sector for the approximation error. Thus, a system analysis and a controller design by SOS techniques are possible. Moreover, we will observe connections to the polynomial approximation approach shown in Section~\ref{Sec_SRNL}.

\subsection{Embedding into linear parameter-varying systems}

In Section~\ref{Sec_LPV}, we will report data-driven system analysis and controller design for a nonlinear system by combining data-driven methods for linear parameter-varying (LPV) systems and embedding nonlinear systems into LPV systems. Thereby, this paradigm provides a data-driven system analysis and controller design by SDPs for nonlinear systems. In contrast to the local approximations by polynomials, LPV systems provide a global linearization of the nonlinear system. However, for the embedding, a known function basis of the scheduling map or the velocity-form of the underlying nonlinear system is required. The LPV representation of the nonlinear system is not tight as the scheduling parameter can change independently of the state and input. 

For LPV systems, \cite{LPV_FundamentalLemma} introduces a fundamental lemma for verifying dissipativity properties \citep{LPV_dissi} and predictive control \citep{LPV_PC}. Moreover, \cite{LPV_direct} provides a representation of open-loop and closed-loop LPV systems from noise-free trajectories. \cite{LPV_set_membership} and \cite{LPV_SetMem2} introduce a set-membership description for LPV systems from noisy data for the controller design with stability and performance guarantees based on SDPs. 

Under the assumption that a function basis of the scheduling map is known, the unknown nonlinear dynamics can be written as an LPV system. Since the nonlinear dynamics is contained within the solution of the LPV system, \cite{LPV_nonlinear} obtains a data-driven LPV controller that stabilizes the underlying nonlinear system. Alternatively, if the function basis of the velocity-form of a nonlinear system is known, then the velocity-form can be embedded into the data-driven framework of LPV systems \citep{LPV_nonlinear2}.

Related to exploiting the extended linearization of a nonlinear system as in \cite{LPV_nonlinear}, the authors \cite{Sznaier} suggest to compute in each time instance a control policy for the frozen system matrices of the nonlinear extended linearization. Together with enforcing a decrease of the Lyapunov function along the closed-loop trajectory, the iterative scheme guarantees closed-loop stability by solving a finite set of LMIs in each time step. The inference on the extended linearization is obtained under a known function basis of the nonlinear dynamics and from data subject to noise.

\subsection{Koopman lifting}

The Koopman operator \citep{Koopman_Survey} provides an exact description of the nonlinear dynamics by a single, but infinite-dimensional, (bi-)linear system. To this end, the time-evolution of the states is observed through the lens of observables. By a finite dictionary of observables, a finite dimensional system representation is obtained from finite data using extended dynamics mode decomposition (EDMD). While the emerging \TM{estimation} error is neglected in many existing results, we will focus in Section~\ref{Sec_Koopman} on works incorporating this error. Thereby, a controller design with guarantees is achieved.

Combined with the data-based validation of the \TM{estimation} error, \cite{Koopman_Guarantees} provides a linear robust predictive control scheme with closed-loop guarantees. \cite{Strasser_Koopman} considers a bilinear lifted system and determines a linear feedback of the lifted states that asymptotically stabilizes the nonlinear system. Moreover, a region of attraction w.r.t. the lifted state is guaranteed. 


Motivated by the bilinear model deduced from the lifting, we will shortly report the data-driven control approaches \cite{Bilinear}, \cite{Bilinear2}, and \cite{Strasser_Bilinear} for bilinear systems with unknown system matrices.

\subsection{Approximate nonlinearity cancellation and feedback linearization}

Whereas the previous approaches derive a suitable representation of the system dynamics itself, nonlinearity cancellation and feedback linearization modify the dynamics by a state feedback to obtain an approximately linear system description. Since the feedback linearization includes an input transformation, the dynamics of the original and the transformed linear system differs. Hence, a data-driven system analysis by feedback linearization or nonlinearity cancellation is intrinsically not possible. For more details, we refer to Section~\ref{Sec_FeedbackLin}.

\cite{NonlinearityCancell} locally asymptotically stabilizes an equilibrium by approximately cancelling the nonlinearity of the closed-loop dynamics. To this end, the data-based closed-loop parametrization for LTI systems \citep{PersisLinear} is obtained for the linear part of the closed loop of a nonlinear system with known function basis. Subsequently, a controller that stabilizes the linear part of the closed-loop dynamics and minimizes the influence of the nonlinearity can be determined by an SDP. 

The global feedback linearization of a flat system serves as the basis for an extension of Willems' fundamental lemma \citep{FeedbackLin2} under a known function basis for the input transformation. \cite{FeedbackLin3} relaxes this assumption by incorporating the error for an arbitrary choice of basis functions into data-driven simulation and output-matching control. Further, \cite{FeedbackLin} proposes a robust predictive control scheme for full-state feedback-linearizable nonlinear systems. However, all these results require nonlinear optimization.

%% file: 03_Polynomial_Approx.tex
\section{Data-driven control by polynomial approximation}\label{Sec_PolyApprox}

A polynomial system representation allows for a system analysis and controller design by SDPs via SOS techniques. Moreover, the well-elaborated error bounds for polynomial interpolation \citep{Hermite} can be leveraged by robust control methods to ensure rigorous guarantees despite approximation of the nonlinear dynamics. 

After an elaboration of data-driven control based on polynomial interpolation techniques, works on polynomial, rational, and Lur'e-type systems will be reported. Lastly, we will review a data-based closed-loop representation, which will also be leveraged in some of the later presented frameworks.

\subsection{Data-driven control by polynomial interpolation}\label{Sec_SRNL}

In the sequel, we present the data-driven set-membership approach for system analysis and controller design from \cite{MartinTP1}, \cite{MartinTP2}, and \cite{MartinCS}. There the feasible system set for nonlinear systems is derived by combining noisy data and error bounds for polynomial interpolation. Since the feasible system set is characterized by polynomial inequalities, the subsequent system analysis and state-feedback design boil down to SOS optimization problems. 

We begin with some background on polynomial interpolation. By Taylor's theorem \citep{TaylorRef}, a $k+1$ times continuously differentiable function $f:\mathbb{R}^{n_x}\rightarrow\mathbb{R}$ can be written as $f(x)=T_k(\omega)[f(x)]+R_k(\omega)[f(x)]$ with the TP of order $k$ at $\omega\in\mathbb{R}^{n_x}$
\begin{align*}
	T_k(\omega)[f(x)]=\sum_{|\alpha|=0}^{k}\frac{1}{\alpha !}\frac{\partial^{|\alpha|}f(\omega)}{\partial x^{\alpha}}\left(x-\omega\right)^\alpha,
\end{align*}
with $\alpha!=\alpha_1!\cdots\alpha_{n_x}!$. Moreover, for all $x\in\mathbb{R}^{n_x}$ there exists a $\nu\in[0,1]$ such that
\begin{align*}
	R_k(\omega)[f(x)]=\sum_{|\alpha|=k+1}\frac{1}{\alpha !}\frac{\partial^{k+1}f(\omega+\nu(x-\omega))}{\partial x^{\alpha}}\left(x-\omega\right)^\alpha.
\end{align*}
Since $\nu$ intrinsically depends on $x$, it summarizes the non-polynomial nonlinearity of $f$. Furthermore, its value is typically unknown. Hence, \cite{MartinTP2} assumes an upper bound on the magnitude of the $(k+1)$-th partial derivatives to avoid the computation of $\nu$. 

\begin{assu}[\cite{MartinTP2}]\label{AssBoundDeri}
\hfill Upper\\ bounds $M_{\alpha}\geq0,\alpha\in\mathbb{N}_0^{n_x},|\alpha|=k+1,$ on the magnitude of each $(k+1)$-th order partial derivative of $f$ are known, i.e.,
\begin{align*}
	\bigg|\bigg|\frac{\partial^{k+1}f(x)}{\partial x^{\alpha}}\bigg|\bigg|_2\leq M_{\alpha},\quad \forall x\in\mathbb{R}^{n_x}. 
\end{align*}
\end{assu}
Under Assumption~\ref{AssBoundDeri}, the Lagrange remainder $R_k(\omega)[f(x)]$ can be bounded by
\begin{equation}\label{PolyRemainderBound}
\begin{aligned}
	(R_k(\omega)[f(x)])^2&\leq \sum_{|\alpha|=k+1}\kappa\frac{M_{\alpha}^2}{\alpha !^2}\left(x-\omega\right)^{2\alpha},
\end{aligned}
\end{equation}
where $\kappa\in\mathbb{N}_0$ is equal to the number of $M_{\alpha}\neq0$. Since the right-hand side of \eqref{PolyRemainderBound} does not depend on $\nu$, we conclude that $f$ is contained in the polynomial characterized sector
\begin{align}\label{Set_Mem_TP}
	f(x)\in\left\{T_k(\omega)[f(x)]+R(x):R^2(x)\leq \sum_{|\alpha|=k+1}\kappa\frac{M_{\alpha}^2}{\alpha !^2}\left(x-\omega\right)^{2\alpha}\right\}.
\end{align}

While the investigation of the approximation error for polynomial interpolation is well-studied, \cite{MartinTP2} proposes their application to infer on the interpolation polynomial from noisy data. For that purpose, let noisy samples $\{y_i,x_i\}_{i=1}^S$ with $y_i=f(x_i)+d_i$ be available. The unknown noise realizations $d_i$ satisfy $d_i^2\leq\epsilon^2$. Furthermore, let the TP 
be written as $z(x)^Ta^*$, where ${a^*}\in\mathbb{R}^{n_z}$ summarizes the unknown coefficients of the interpolation polynomial and $z(x)\in\mathbb{R}[x]^{n_z}$ is a vector of linear independent polynomials building a basis for all polynomials with degree less than or equal to $k$. For instance, $z$ could contain all monomials up to degree $k$. 
Following a set-membership procedure \citep{Milanese_SysId}, we characterize all coefficients $a\in\mathbb{R}^{n_z}$ admissible with the data
\begin{equation}\label{Set_mem_coeff}
\begin{aligned}
	&\{a:\exists  {d}_i\in\mathbb{R},i=1,\dots,S, \text{\ satisfying\ } {d}_i^2\leq\epsilon^2  \text{\ and\ }\\
	 &\hspace{3cm}{y}_i=z({x}_i)^Ta+R_k(\omega)[f(x_i)]+{d}_i\}.
\end{aligned}
\end{equation}
Note that \eqref{Set_mem_coeff} includes the coefficients of the true interpolation polynomial ${a^*}$. Furthermore, since the remainder evaluated at the data are unknown, \cite{MartinTP2} (Lemma 2) calculates a superset of \eqref{Set_mem_coeff} based on the error bounds for TPs \eqref{PolyRemainderBound} 
, denoted by $\bar{R}^{\text{poly}}[f(x_i)]$,
\begin{align*}
 	{\Sigma}_{a}= \left\{a:
 	({y}_{i}-z({x}_i)^Ta)^2\leq q({y}_i,{x}_i) ,i=1,\dots,S\right\}
\end{align*}
with $q({y}_i,{x}_i)=\bar{R}^{\text{poly}}[f({x}_i)]+\epsilon^2+2\epsilon\sqrt{ \bar{R}^{\text{poly}}[f({x}_i)]}$. Finally, combining the polynomial sector bound for TPs \eqref{Set_Mem_TP} 
with ${\Sigma}_{a}$ leads to a set membership for the unknown nonlinear function 
\begin{equation}\label{Set_mem_TP}
	f(x)\in\left\{z(x)^Ta+R(x):a\in{\Sigma}_{a},R(x)^2\leq\bar{R}^{\text{poly}}[f({x})]\right\}.
\end{equation}
By the latter, we obtain for the nonlinear function $f$ a representation that is based on noisy data, is polynomial, and does not require a function basis of $f$.\\

Before we continue with the application of the feasible system set \eqref{Set_mem_TP} for data-driven system analysis and control, some comments are appropriate. By Assumption~\ref{AssBoundDeri}, the rate of variation of the nonlinear function $f$ is bounded. Therefore, we can infer the behavior of $f$ in the neighbourhood of a sample. 
A similar causality of the behavior of $f$ and samples is, e.g., also considered in set-membership identification \citep{Milanese_SysId} by a known Lipschitz constant. Since the bounds $M_{\alpha}$ are usually not known, we refer to \cite{MartinTP1} for their estimation by a validation procedure and the applicability for data from a real experiment.

Moreover, \cite{MartinTP1} discusses the choice of $\omega$ and a combining of multiple local polynomial interpolations, i.e., a piecewise polynomial approximation, to reduce the approximation error. To reduce the computation burden, \cite{MartinTP2} (Proposition 1) suggests to compute the ellipsoidal outer approximation
\begin{equation}\label{Quad_const}
	\left\{a:\begin{bmatrix}I\\a^T\end{bmatrix}^TQ\begin{bmatrix}I\\a^T\end{bmatrix}\preceq0 \right\}
\end{equation}
of ${\Sigma}_{a}$. Furthermore, \cite{MartinCS} investigates a Frequentist and Bayesian treatment for Gaussian distributed $d_i$, following the lines of \cite{Umenberger}.

Next, we show how to verify system properties for a general nonlinear system by the derived set membership \eqref{Set_mem_TP}. For that purpose, \cite{MartinTP2} analyzes an unknown nonlinear continuous-time system 
\begin{equation*}
	\dot{x}(t)=f(x(t),u(t))
\end{equation*}
with $k+1$ times continuously differentiable function $f=\begin{bmatrix}f_1 & \cdots & f_{n_x}\end{bmatrix}:\mathbb{R}^{n_x+n_u}\rightarrow\mathbb{R}^{n_x}$ and noisy data $\{\dot{x}_i,x_i,u_i\}_{i=1}^S$ with $\dot{x}_i=f(x_i,u_i)+d_i$. Note that $d_i$ might include estimation errors of the time-derivatives of the states. Applying a polynomial interpolation for each element of $f$ yields 
\begin{align*}
	f(x,u)=&\begin{bmatrix}z_1(x,u)^Ta_1^*+R[f_1(x,u)]\\\vdots\\ z_{n_x}(x,u)^Ta_{n_x}^*+R[f_{n_x}(x,u)]\end{bmatrix}\\
	=&\underbrace{\begin{bmatrix}z_1(x,u)^TS_1\\\vdots\\ z_{n_x}(x,u)^TS_{n_x}\end{bmatrix}}_{=:Z(x,u)}a^*+\underbrace{\begin{bmatrix}R[f_1(x,u)]\\\vdots\\ R[f_{n_x}(x,u)]\end{bmatrix}}_{=:R[f(x,u)]},
\end{align*}
with $a_i^*=S_ia^*$, and $R[f(x,u)]^TR[f(x,u)]\leq\sum_{i=1}^{n_x}\bar{R}^{\text{poly}}[f_i({x,u})]$. 
As shown in \cite{MartinTP1} (Remark 2) and \cite{MartinCS}, summarizing all unknown coefficients of the polynomial interpolations into $a^*$ and considering $z_i$ for each row of $f$ can be leveraged to include interpolation polynomials of different orders and prior knowledge on the structure of $f$. 

For the available data $\{\dot{x}_i,x_i,u_i\}_{i=1}^S$, we pursue the described derivation of ${\Sigma}_{a}$ to conclude that $f(x,u)$ is contained in the feasible system set  
\begin{equation}\label{Set_mem_SysDyn}
\begin{aligned}
	&\left\{Z(x,u)^Ta+R(x,u):a\in{\Sigma}_{a},\phantom{\sum_{i=1}^{n_x}}\right.\\&\hspace{2cm} \left.R(x,u)^TR(x,u)\leq\sum_{i=1}^{n_x}\bar{R}^{\text{poly}}[f_i({x,u})]\right\}.
\end{aligned}	
\end{equation}
Given the quadratic description \eqref{Quad_const} for ${\Sigma}_{a}$, the feasible system set can also be written as an LFR \citep{SchererLMI}
\begin{equation*}
\begin{aligned}
\begin{bmatrix}
	\dot{x}(t)\\q(t)
\end{bmatrix}&=\begin{bmatrix}0 & 0 & I &I\\ \begin{bmatrix}I\\0\end{bmatrix} & \begin{bmatrix}0\\I\end{bmatrix} & 0 &0
\end{bmatrix}\begin{bmatrix}
x(t)\\ u(t)\\ w_1(t)\\ w_2(t)
\end{bmatrix},\\
w_1(t) &= Z(q(t))^Ta, w_2(t) = R(q(t)), 
\end{aligned}
\end{equation*}
with $w_2(t)^Tw_2(t)\leq\sum_{i=1}^{n_x}\bar{R}^{\text{poly}}[f_i({q(t)})]$ and \begin{equation*}
	\begin{bmatrix}S_i^Tz_i(q(t))\\w_{1,i}(t)\end{bmatrix}^TQ\begin{bmatrix}S_i^Tz_i(q(t))\\w_{1,i}(t)\end{bmatrix}\leq0,
\end{equation*}
for $w_1=\begin{bmatrix}w_{1,1} & \cdots&w_{1,n_x}\end{bmatrix}^T$. 

Although the uncertainty descriptions of $w_1$ and $w_2$ are polynomial in $q$, we can apply the LMI-based robust control framework of \cite{SchererLMI} to verify system properties for all systems within the set membership \eqref{Set_mem_SysDyn}. If all systems satisfy a certain property, then also the ground-truth system fulfils the property as it is contained within \eqref{Set_mem_SysDyn}. Due to the polynomial bounds on the uncertainties $w_1$ and $w_2$, the verification boils down to an SOS condition. By an additional S-procedure, \cite{MartinTP2} and \cite{MartinTP1} verify dissipativity and incremental dissipativity properties, respectively, for all trajectories of the unknown nonlinear system staying within a compact set. A regional analysis is meaningful due to the non-global inferences on the dynamics from data. 
Moreover, the results from \cite{MartinIQC} (arXiv:2103.10306v3, Section 5.B) can be applied here to determine optimal IQCs to gather tighter system properties than by simple dissipativity properties. 

\cite{MartinCS} elaborates the set membership \eqref{Set_mem_SysDyn} with a TP at $\omega=0$ for each row of an input-affine nonlinear system. Following the robust control framework of \cite{SchererLMI}, Theorem 8 of \cite{MartinCS} formulates the conditions for data-driven dissipativity verification in the dual space. Thereby, a state-feedback design by SOS optimization with quadratic performance guarantees is possible. In contrast to the TP approach from \cite{DePersis2}, the authors \cite{MartinCS} incorporate the remainder of polynomial interpolation into the controller synthesis to obtain globally asymptotically stabilizing controllers.

\begin{figure}
	\begin{center}
		\input{Img/Indirect_scheme}
	\end{center}
	\caption{Graphical illustration of data-driven system analysis and control using polynomial interpolation.}
	\label{Fig.Scheme}
\end{figure}
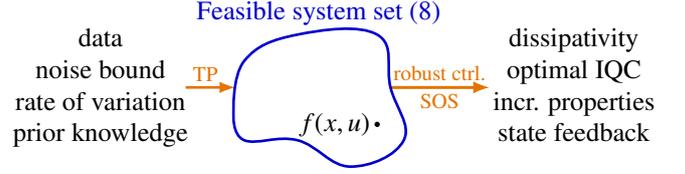

Figure~\ref{Fig.Scheme} summarizes the data-driven system analysis and control by polynomial interpolation.
While this framework seems to be flexible, we also identify some open questions and problems. Among others, the impact of the chosen polynomial basis for $z_1,\dots,z_{n_x}$ should be analysed. Here, the basis should be optimized for the controller design. Hence, the basis is ideally optimized during the controller synthesis. 
Since the the presented procedure can also be applied for more general polynomial interpolations, e.g., Hermite polynomials \citep{Hermite}, this investigation might reduce the conservatism by means of more regional polynomial approximations.

\subsection{Polynomial, rational, and Lur'e-type systems}\label{Sec_Poly_rational}

\begin{figure}
	\centering
	\begin{tikzpicture}
	\node[inner sep=0pt] (p1) at (-5.1,0) {};
	\node[inner sep=0pt] (p1) at (0,0) {\includegraphics[width=.5\textwidth]{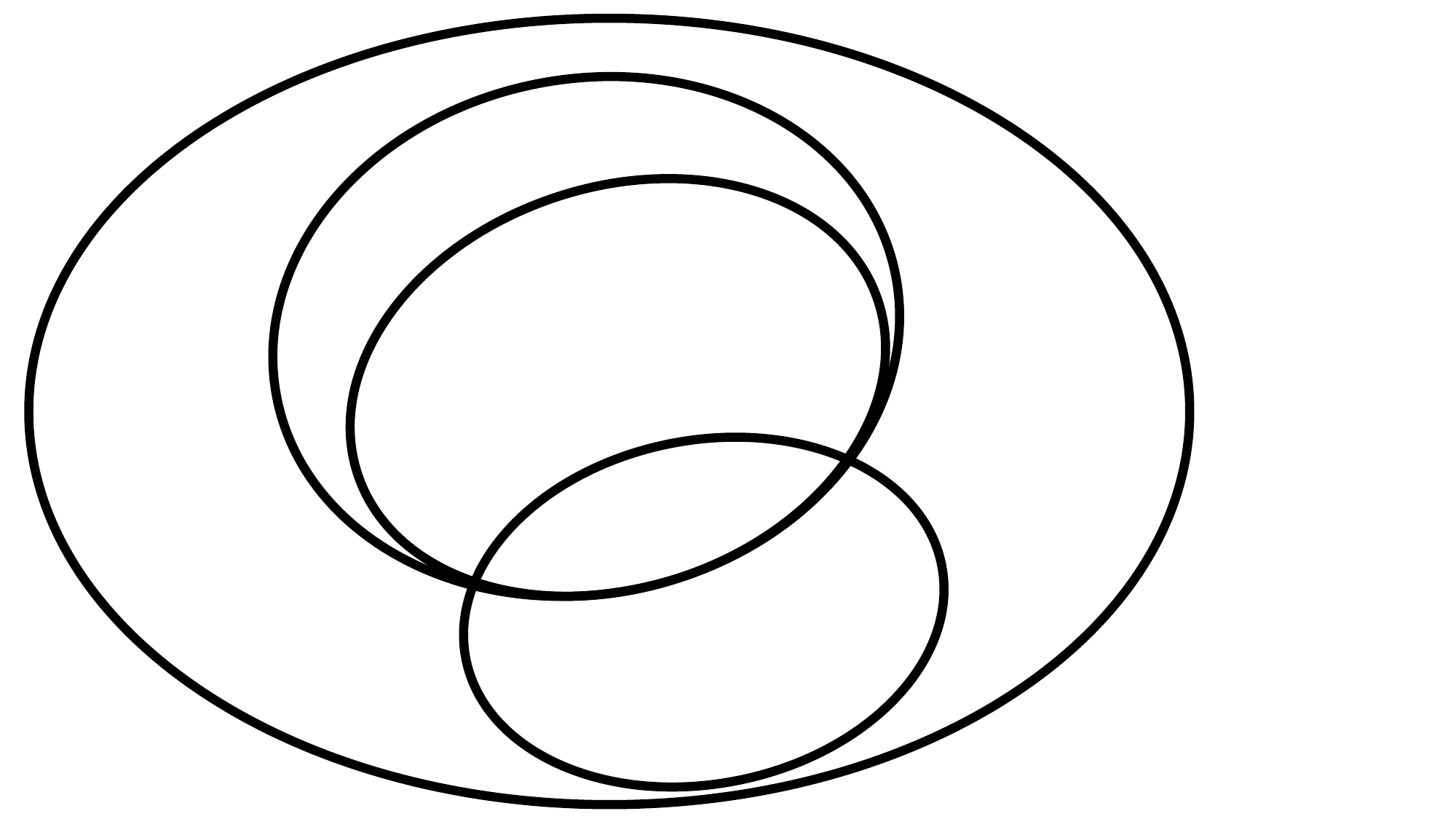}};
	
	\node[inner sep=0pt] (p1) at (-0.5,-0.7) {LTI};
	\node[inner sep=0pt] (p1) at (-1.1,1.7) {Rational};
	\node[inner sep=0pt] (p1) at (-0.7,0.3) {Polynomial};
	
	\node[inner sep=0pt] (p1) at (-0.5,-1.7) {Lur'e};
	\node[inner sep=0pt] (p1) at (1.95,0.0) {\begin{tabular}{c}
		Polynomial\\interpolation\end{tabular}
		 };

	\end{tikzpicture}
	\caption{Relation of polynomial interpolation approach from Section~\ref{Sec_SRNL} and approaches for subclasses from Section~\ref{Sec_Poly_rational}.}
	\label{Fig.PolyInter_subClass}
\end{figure}
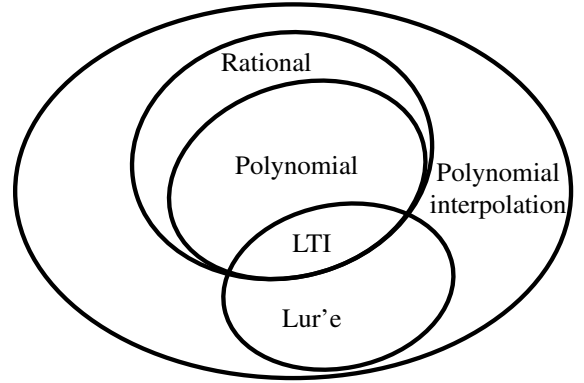

We summarize in the following the data-driven set-membership approaches for polynomial, rational, and Lur'e-type systems. These classes of nonlinear systems are strongly related to the system representation by polynomial interpolation as shown in Figure~\ref{Fig.PolyInter_subClass}. 
Polynomial systems comprises, for example, fluid dynamics \citep{Polynomial_sys} and robotics \citep{Poly_Robotic}, whereas rational dynamics, e.g., biochemical reactors \citep{Rational_sys}. Lur'e systems \citep{Khalil} (Chapter 10.1) are LTI systems including a sector bounded nonlinearity, and thus can comprise exemplarily Lipschitz bounded nonlinearities and recurrent neural networks \citep{Persis_quad_constr}.

The literature on data-driven system analysis and control for polynomial systems mostly studies a continuous-time system with polynomial dynamics $\dot{x}(t)=Az(x(t),u(t))$ or input-affine polynomial dynamics $\dot{x}(t)=A\bar{z}(x(t))+BW(x(t))u(t)$, respectively. While the coefficient matrices $A\in\mathbb{R}^{n_x\times n_z}$ and $B\in\mathbb{R}^{n_x\times n_B}$ are unknown, a set of data $\{\dot{x}_i,x_i,u_i\}_{i=1}^S$ is available and a vector $z\in\mathbb{R}[x,u]^{n_z}$ or a vector $\bar{z}\in\mathbb{R}[x]^{n_z}$ and a polynomial matrix $W\in\mathbb{R}[x]^{n_B\times n_u}$, respectively, are known. The latter can be satisfied from knowledge of an upper bound on the degree of the polynomial dynamics.

Clearly, considering polynomial systems generalizes the results for LTI systems \citep{vanWaarde,LFT_SA}. On the other hand,
the polynomial problem setup is included in the data-driven framework of polynomial interpolation in Section~\ref{Sec_SRNL}: For vanishing approximation error $\bar{R}^{\text{poly}}[f({x,u})]=0$, we can proceed as in Section~\ref{Sec_SRNL} to derive a set membership for the unidentified coefficients and, subsequently, to verify system properties and to design a state feedback. Hence, we directly report the examined data-driven control problems for polynomial systems.

Based on SOS optimization, {\cite{MartinPoly}} and {\cite{MartinIQC}} investigate the verification of dissipativity and more general time domain hard IQCs with optimized linear filter, respectively. Whereas most model-based and data-driven results, which are based on SOS optimization, consider continuous-time systems, \cite{MartinPoly} and \cite{MartinIQC} investigate discrete-time polynomial systems. Thereby, the data of the time-derivative of the states are avoided. This advantage comes at the cost of being restricted to quadratic storage functions \citep{MartinIQC} or polynomials of higher degree due to the function decomposition $V(f(x,u))$ with a polynomial storage function $V(x)$. Moreover, polynomial discrete-time systems often require a regional analysis because a time discretization of a globally asymptotically stable system can lead to a locally asymptotically stable discrete-time system. For instance, the Euler time discretization with time step $T>0$ of $\dot{x}=-x^3$ yields $x(t+1)=x(t)-Tx(t)^3$. However, the state of the discretized system 
	tends to infinity for any initial condition $||x(0)||_2>\sqrt{2/T}$.


Furthermore, data-driven controller design for global stabilization for continuous-time polynomial systems is examined: {\cite{DePersis1}} derives a polynomial state feedback from data solving a single SOS condition. To this end, an energy-bounded noise of the data immediately leads to a set membership for the unknown coefficient matrices, which can be exploited by means of the non-conservative S-lemma from \cite{vanWaarde}. The same control problem for pointwise-in-time-bounded noise is solved in {\cite{Petersen}} by the non-conservative Petersen's lemma and an ellipsoidal outer approximation of the set of coefficients consistent with the data. However, the controller design requires alternating SOS optimization due to the bilinearity of optimization variables. {\cite{FarkasLem}} proposes a rational state feedback by relating the feasible system set with the set of systems with converging trajectories via Farkas' lemma. Therefore, the weaker stability of Rantzer's dual Lyapunov theory \citep{Ranther_stab} is considered. Farkas' lemma leads to a condition that can be relaxed to an SDP by SOS relaxation and the nuclear norm relaxation for rank conditions. Summarized, the literature on data-driven state-feedback design considers a variety of robust control techniques to stabilize all systems consistent with the data. However, we assess the robust control framework by LFRs \citep{SchererLMI} to constitute one of the most appealing ones. In particular, it can also be applied for polynomial problems and leads to SDPs despite multiple uncertainty channels and performance criteria. 

In contrast to global stabilization, {\cite{Poly_Invariance}} and {\cite{SafeControl2}} investigate the control of unknown polynomial systems including safety conditions rather than stability conditions. To this end, \cite{Poly_Invariance} optimizes over a polynomial state-feedback law to enlarge the size of an invariant set. 
The feasible system set is incorporated into the invariance condition by Young's relation \citep{Young}. But in fact, as in the model-based case, the derived invariance conditions are bilinear, and thus have to be solved by alternating SOS optimization. To circumvent the nonconvexity, \cite{SafeControl2} proposes a density function formulation based on the dual Lyapunov method \citep{Ranther_stab} and to proceed along the lines of \cite{FarkasLem}. The obtained rational state feedback keeps the closed-loop trajectories of all polynomial systems admissible with the data outside of an unsafe set.\\

Instead of polynomial systems, {\cite{Strasser}} presents a data-driven feedback design for continuous-time rational systems.  Multiplication of the rational dynamics by all denominators of the dynamics yields a problem formulation akin to that in the polynomial case. However, the problem emerges that the additive noise does not affect the data through the original rational dynamics but the polynomial reformulation. 
The remaining procedure follows a set-membership approach combined with robust control techniques from \cite{SchererLMI}. Thereby, a state feedback is designed with stability and performance guarantees by solving an SOS problem.\\

Related to a linear TP with bounded Lagrange remainder, the Lur'e problem considers an LTI system with a sector bounded nonlinearity. In a data-driven context, {\cite{Berberich}} proposes a flexible multiplier LMI-framework to combine data, prior information on the LTI system, and bounded and measurable nonlinearities. To obtain a feasible system set for the LTI system with sector bounded and measurable nonlinearity, {\cite{Findeisen}} solves an optimal control problem from MPC with state and input constraints by application of the S-lemma from \cite{vanWaarde}. {\cite{Input_constr}} solves a similar problem by Young's relation.

\subsection{A data-driven closed-loop characterization}\label{Sec_Direct}

\cite{PersisLinear} examines besides the data-driven stabilization of LTI systems also the nonlinear case. In contrast to the set-membership approaches of Section~\ref{Sec_SRNL} and \ref{Sec_Poly_rational}, \cite{PersisLinear} proposes to directly characterize the closed loop by data matrices. This result has also inspired further works reported, among others, in this section. 

\cite{PersisLinear} (Section 5.2) considers the stabilization of an unknown nonlinear discrete-time system $x(t+1)=f(x(t),u(t))$. To apply their results from LTI systems, the authors propose a Taylor linearization around the known equilibrium point $(x_\text{e},u_\text{e})=(0,0)$. This leads to the linearized system dynamics
\begin{equation}\label{TP_linearized}
	x(t+1)=Ax(t)+Bu(t)+d(t),
\end{equation}    
with $A=\frac{\partial f(0,0)}{\partial x}$, $B=\frac{\partial f(0,0)}{\partial u}$, and remainder $d$. Clearly, \eqref{TP_linearized} corresponds to the TP approximation from Section~\ref{Sec_SRNL} for $k=1$ and $\omega=0$. For data $\{x_i,u_i\}_{i=1}^S$ satisfying $x_{i+1}=Ax_i+Bu_i+d_i$, \cite{PersisLinear} defines the data-dependent matrices $X=\begin{bmatrix}x_1 & \cdots&x_{S-1}\end{bmatrix}$, $U=\begin{bmatrix}u_1 & \cdots&u_{S-1}\end{bmatrix}$, $X^+=\begin{bmatrix}x_2 & \cdots&x_{S}\end{bmatrix}$, and $D=\begin{bmatrix}d_1 & \cdots&d_{S-1}\end{bmatrix}$.

\begin{assu}[\cite{PersisLinear} (Assumption 5)]\label{AssBoundLinearization}
Let a constant $M>0$ with $DD^T\preceq MX^+{X^+}^T$ be known.
\end{assu}
Assumption~\ref{AssBoundLinearization} cumulatively bounds the whole sequence of the remainder $d_1,\dots,d_{S-1}$ by a single constraint, rather than for each realization $d_i$ separately as in Section~\ref{Sec_SRNL}. Furthermore, in contrast to the quadratically increasing bound of the Lagrange remainder \eqref{PolyRemainderBound}, Assumption~\ref{AssBoundLinearization} supposes a bound on the remainder that increases linearly. Nevertheless, $M$ of Assumption~\ref{AssBoundLinearization} corresponds to similar insights as $M_\alpha$ in Assumption~\ref{AssBoundDeri}.

For deriving a stabilizing linear state feedback $u=Kx$, observe that the data matrices satisfy $X^+=AX+BU+D$. Hence, the linear part of the closed-loop dynamics can be characterized based on the data matrices
\begin{align}\label{DD_CL}
	A+BK=\begin{bmatrix}A & B\end{bmatrix}\begin{bmatrix}I\\K\end{bmatrix}=\begin{bmatrix}A & B\end{bmatrix}\begin{bmatrix}X\\U\end{bmatrix}G=(X^+-D)G,
\end{align}
with $G$ satisfying
\begin{equation}\label{consitency_cond}
	\begin{bmatrix}I\\K\end{bmatrix}=\begin{bmatrix}X\\U\end{bmatrix}G.
\end{equation}
The persistence of excitation (PE) condition that $\begin{bmatrix}X\\U\end{bmatrix}$ has full row rank \citep{PersisLinear} (Assumption 4) guarantees that $G$ always exists. Note that a data-based description of the closed-loop dynamics can also be retrieved by the relationship $\begin{bmatrix}A & B\end{bmatrix}=(X^+-D)\begin{bmatrix}X\\U\end{bmatrix}^\dagger$, which corresponds to $G=\begin{bmatrix}X\\U\end{bmatrix}^\dagger\begin{bmatrix}I\\K\end{bmatrix}$. Instead of one specific $G$, the closed-loop description \eqref{DD_CL} with consistency condition \eqref{consitency_cond} provides for a fixed $K$ a whole subspace in terms of $G$. Hence, optimizing over $G$ typically results in non-unique solutions. As shown in \cite{Regularization}, it is advantageous to regularize $G$ to single out solutions that are robust regarding noise.

Along Theorem 6 of \cite{PersisLinear}, the origin of the nonlinear system is asymptotically stable if the linear part of the closed loop \eqref{DD_CL} is asymptotically stable, i.e., if there exists a Lyapunov matrix $P\succ0$ and matrix $G$ such that
\begin{equation*}
	(A+BK)P(A+BK)^T-P=((X^+-D)G)P((X^+-D)G)^T-P\prec0.
\end{equation*}
Since the remainder matrix $D$ is unknown, Theorem 5 of \cite{PersisLinear} uses Young's relation \citep{Young} to ensure by an SDP that this stability condition holds for all $D$ satisfying Assumption~\ref{AssBoundLinearization}. By optimizing over $P$ and $G$, the state-feedback matrix $K$ can be recovered from \eqref{consitency_cond}. 

\cite{PersisLinear} describes the uncertainty from the unknown remainder directly by the matrix $D$. Thereby, \cite{PersisLinear} obtains a parametrization $G$ of the to-be-optimized controller that increases with the number of data. In contrast, the set-membership approach in Section~\ref{Sec_SRNL} translates the uncertain remainder into an uncertainty of the coefficients of the TP. This leads to a controller design that does not scale with the number of samples. Moreover, incorporating the remainder of the TP approximation into the controller synthesis \citep{MartinCS} enables a global stabilization. Throughout this article, we will see further results based on a set-membership procedure or the closed-loop characterization from \cite{PersisLinear} exhibiting similar properties.\\   

Some of the data-based results for nonlinear systems inspired by \cite{PersisLinear} are mentioned next. \cite{DePersis1} includes an extension for polynomial systems. \cite{ExpDesign} proposes an experiment design by scaled input sequences to ensure Assumption~\ref{AssBoundLinearization} and the PE condition for the data from a nonlinear system. Furthermore, \cite{DDC_orbit} stabilizes a periodic orbit using a Pyragas-type control law and a system representation for the periodic time-evolution of the states. While the to-be-stabilized orbit becomes an equilibrium of the periodic system description, \cite{DDC_orbit} also solves the problem when the orbit is not precisely known. \cite{Persis_quad_constr} stabilizes Lur'e-type systems $x(t+1)=Ax(t)+B(u(t))+f(t,x(t))$ using SDPs with LMI constraints. To this end, a quadratic constraint and samples of the nonlinearity $f$ need to be known, analogously to the approaches for Lur'e systems mentioned in Section~\ref{Sec_Poly_rational}.\\ 

Although not directly related to the presented closed-loop representation of \cite{PersisLinear}, we would also like to mention \cite{BerberichNL}. The authors also use a data-based inference on the behavior of the affine Taylor linearization, but in the context of predictive control. To this end, online updated data is combined with the fundamental lemma of \cite{Willems} to obtain a data-driven system parametrization of the nonlinear system in the neighbourhood of the current state. 
Thereby, the MPC scheme boils down to solving a convex quadratic program in each time instance.

%% file: Img/Indirect_scheme.tex
\begin{tikzpicture}[scale=0.7]

\draw[fill] (3.4,0.8) circle [radius=0.04];
\path (2.6,0.8) node[](S2) {\textcolor{black}{$f(x,u)$}}; 	
\path (2.3,2.9) node[](S2) {\textcolor{blue!80!black}{Feasible system set (8)}};

\path (-1.8,1.5) node[](S3) {\begin{tabular}{c}
	\textcolor{black}{data}\\\textcolor{black}{noise bound}\\\textcolor{black}{rate of variation}\\\textcolor{black}{prior knowledge}\end{tabular}};
\path (7.1,1.5) node[](S1) {\begin{tabular}{c}	\textcolor{black}{dissipativity}\\\textcolor{black}{optimal IQC}\\\textcolor{black}{incr. properties}\\\textcolor{black}{state feedback}\end{tabular}};

\draw[latex-,line width=1pt,orange!90!black] (5.5,1.5) -- (3.65,1.5)node[pos=0.5, below=2pt, inner sep=0pt](HP4){\footnotesize\textcolor{orange!90!black}{SOS}};
\draw[latex-,line width=1pt,orange!90!black] (5.5,1.5) -- (3.65,1.5)node[pos=0.5, above=2pt, inner sep=0pt](HP4){\footnotesize\textcolor{orange!90!black}{robust ctrl.}};

\draw[latex-,line width=1pt,orange!90!black] (0.73,1.5) -- (-0.2,1.5)node[pos=0.6, above=2pt, inner sep=0pt](HP4){\footnotesize\textcolor{orange!90!black}{TP}};

\draw[draw=blue!80!black,line width=1pt]  plot[smooth, tension=.7] coordinates {(0.7,1) (0.9,2) (1.9,2.5) (3.4,2.5) (3.65,1.5) (3.9,0.5) (2.9,0) (1.9,0.5) (0.9,0.5) (0.7,1)};

\end{tikzpicture}

%% file: 04_GP.tex
\section{GP and kernel ridge regression for data-driven control}\label{Sec_GP}

GP and kernel ridge regression constitute well-established techniques in machine learning to approximate a nonlinear function from data and to predict its outputs for unseen inputs \citep{GP_ML} and \citep{Bishop} (Section~6). Both regression methods are equipped with an uncertainty measure given by an upper bound for their approximation error. However, the regression solution and its error bound are usually nonlinear due to the nonlinear kernel functions. For that reason, the research direction, reported in this section, tackles the challenge to identify a representation of the regression that is suitable for system analysis and controller synthesis using SDPs. Thereby, these data-based methods can be leveraged for a wide range of control problems compared to the specific control schemes, e.g., using GP models with feedback linearization \citep{Feedbacklin_alternative1} and \citep{Feedbacklin_alternative2}.

\subsection{Kernel ridge regression}

For kernel ridge regression, we first introduce the notion of kernel functions and their reproducing kernel Hilbert space (RKHS) to efficiently solve the following infinite-dimensional regression problem: For an unknown nonlinear function $f:\mathbb{X}\subseteq\mathbb{R}^{n_x}\rightarrow\mathbb{R}$, let the data points $\{y_i,x_i\}_{i=1}^{S}$ with $y_i=f(x_i)$ be available. Then, we want to find the solution of
\begin{equation}\label{Kernel_regression}
	\min_{\mu\in\mathcal{H}}\sum_{i=1}^{S}(y_i-\mu(x_i))^2+\lambda\,||\mu||_\mathcal{H}^2,
\end{equation}
where $\lambda>0$ is a regularization parameter, $\mathcal{H}$ a real Hilbert space of functions $\mu:\mathbb{X}\rightarrow\mathbb{R}$, and $||\cdot||_\mathcal{H}$ the associated function norm, i.e., $||\mu||_\mathcal{H}=\sqrt{\langle \mu,\mu\rangle_\mathcal{H}}$.

To solve this least-squares-error (LSE) problem, let a function $k:\mathbb{X}\times\mathbb{X}\rightarrow\mathbb{R}$ be a kernel if it is symmetric $k(x,x')=k(x',x)$ for all $x,x'\in\mathbb{X}$ and positive definite, i.e., $\sum_{i,j=1}^{n}\alpha_i\alpha_j k(x_i,x_j)\geq0$ for all $n\in\mathbb{N}$, $x_1,\dots,x_n\in\mathbb{X}$, and $\alpha_1,\dots,\alpha_n\in\mathbb{R}$ \citep{Steinwart} (Theorem 4.16). Kernels naturally emerge in finite dimensional linear regression problems and allow a computationally efficient evaluation. Indeed, the solution of a linear regression includes the evaluation of the scalar product of the feature mapping $Z(x):\mathbb{X}\rightarrow\mathbb{R}^{n_f}$, with often large $n_f$. At the same time, there exists for many $Z(x)$ a kernel $k$ such that $k(x,x')=Z(x)^TZ(x')$ \citep{Steinwart} (Definition 4.1). 

The special class of reproducing kernels and their associated RKHS will be crucial to solve \eqref{Kernel_regression}. According to \cite{Steinwart} (Def. 4.18), a kernel $k:\mathbb{X}\times\mathbb{X}\rightarrow\mathbb{R}$ is a reproducing kernel of a real Hilbert space $\mathcal{H}$ of functions $f:\mathbb{X}\rightarrow\mathbb{R}$ if $k(x,\cdot)\in\mathcal{H}$ and $f(x)=\langle f,k(x,\cdot)\rangle_\mathcal{H}$ for all $x\in\mathbb{X}$ and $f\in\mathcal{H}$. The Hilbert space $\mathcal{H}$ is called RKHS of kernel $k$. 
 
In the sequel, let a kernel $k$ and its RKHS $\mathcal{H}$ be given and the following assumption be satisfied. 
\begin{assu}[\cite{Kernel_assum}]\label{AssKernel}
	Let $f\in\mathcal{H}$ and let an upper bound $M$ on $||f||_{\mathcal{H}}$ be known. 
\end{assu}
Under Assumption~\ref{AssKernel}, the regression problem \eqref{Kernel_regression} is called kernel ridge regression and the representer theorem \citep{RepresenterTheorem} (Theorem 3.4) provides its explicit solution
\begin{equation}\label{Kernel_apporx}
	\mu(x)=y_X(\lambda I+K_X)^{-1}K(x)\in\mathcal{H},
\end{equation} 
with $y_X=\begin{bmatrix}y_1 & \cdots & y_S\end{bmatrix}$, $K(x)=\begin{bmatrix}k(x,x_1)&\cdots & k(x,x_S)\end{bmatrix}^T$, and Gram matrix $K_X$, i.e., its $(i,j)$-th element corresponds to $k(x_i,x_j)$. Furthermore, \cite{Persis_kernel} derives the following approximation error
\begin{equation*}
	||f(x)-\mu(x)||_2\leq M\sqrt{k(x,x)-K(x)^T\hat{K}_X^{-1}K(x)},
\end{equation*}
for all $x\in\mathbb{X}$ and with $\hat{K}_X=(\lambda I+K_X)(2\lambda I+K_X)^{-1}(\lambda I+K_X)$. Thus, we obtain the set membership
\begin{equation}\label{Set_mem_kernel}	
\begin{aligned}
	f(x)\in\Sigma_{\text{ker}}=\left\{\mu(x)+R(x): \mu(x)=y_X(\lambda I+K_X)^{-1}K(x)\right.\\ \left.\text{ and } R(x)^2\leq M^2(k(x,x)-K(x)^T\hat{K}_X^{-1}K(x)) \right\}.
\end{aligned}
\end{equation}

By \eqref{Set_mem_kernel}, we can conclude that the unknown function $f$ is contained within a sector with centre $\mu(x)$ and width $M\sqrt{k(x,x)-K(x)^T\hat{K}_X^{-1}K(x)}$. In contrast to the set membership \eqref{Set_mem_TP} from polynomial approximation, the description of $\Sigma_{\text{ker}}$ is usually nonlinear as the kernel $k$ is typically nonlinear, for example, compare the list of kernels in Section~2.1 of \cite{RepresenterTheorem}. Hence, a direct application of $\Sigma_{\text{ker}}$ for a system analysis or a controller design by SDPs is not possible. Furthermore, Assumption~\ref{AssKernel} requires a bound on the norm associated to the kernel of the nonlinear function, while Assumption~\ref{AssBoundDeri} a bound on high order partial derivatives. Both insights enable one to bound the difference between the nonlinear function and the approximation by the kernel regression or the polynomial interpolation. With the approximation error for $f\notin\mathcal{H}$ from \cite{Fiedler2}, both approaches do not require knowledge of a function basis of the underlying nonlinear function. A data-driven inference on Assumption~\ref{AssKernel} is examined in \cite{Bound_kernel}. Moreover, an approximation error for noisy data is investigated in \cite{Kernel_assum}. Finally, since the kernel and the functions of its RKHS can share desired properties \citep{Prior_kernel}, prior knowledge on $f$ can be included to reduce conservatism and improve the data efficiency.

\subsection{GP regression}

GP regression is another non-parametric regression method to infer an unknown nonlinear function by refining a prior belief by a set of noisy data. Hence, GP regression is a Bayesian method, where the underlying data generation is not drawn from a fixed function (Frequentist statistics) but from a stochastic process $\mathcal{F}$, i.e., a distribution over functions. 

In the sequel, let a random vector $X$ that is Gaussian distributed with mean $\mu$ and covariance matrix $\Xi\succ0$ be denoted by $X\sim\mathcal{N}(\mu,\Xi)$. A GP is a collection of random variables such that any finite subset is Gaussian distributed \citep{GP_ML}. This specific stochastic process is uniquely defined by its mean function $m:\mathbb{X}\rightarrow\mathbb{R}$ and its covariance kernel $k(x,x')$. If the prior distribution $\mathcal{F}$ is a GP with zero mean and covariance $k$ and the data $\{y_i,x_i\}_{i=1}^S$ with $y_i=\mathcal{F}(x_i)+d_i$ and independently identically distributed (iid) $d_i\sim\mathcal{N}(0,\sigma^2)$ is available, then the posterior $\mathcal{F}_\text{post}$ is again a GP with mean 
\begin{equation}\label{post_GP_mean}
	\mu_\text{post}(x)=y_X(\sigma^2 I+K_X)^{-1}K(x)
\end{equation}
and covariance
\begin{equation}\label{post_GP_cov}
	k_\text{post}(x,x')=k(x,x')-y_X(\sigma^2 I+K_X)^{-1}K(x)
\end{equation}
\citep{GP_ML}. Here $y_X, K_X$, and $K(x)$ are defined as for the kernel ridge regression. Thus, the kernel ridge regression solution \eqref{Kernel_apporx} is the same as the posterior mean. Intuitively, for large regularization parameters $\lambda$ or large noise variance $\sigma^2$, the regression loosely fits the samples or does not trust the data, respectively. For more details on the connections between GP and kernel ridge regression, we refer to \cite{RepresenterTheorem}. Moreover, the posterior variance $k_\text{post}(x,x)$ measures the uncertainty of the inference on the data-generating distribution $\mathcal{F}$. 

In the context of robust control, the ground truth dynamics might not be a stochastic process $\mathcal{F}$ but a function $f:\mathbb{X}\rightarrow\mathbb{R}$. Since the posterior mean and covariance are obtained from Bayesian methods, a Frequentist bound on $f$ can not be derived directly from the previous Bayesian treatment. To infer a Frequentist bound, let data $\{y_i,x_i\}_{i=1}^S$ with $y_i=f(x_i)+d_i$ and iid $d_i\sim\mathcal{N}(0,\sigma^2)$ be given. Then we can rely on Assumption~\ref{AssKernel} to compute a bound on the approximation error of the form
\begin{equation}\label{appro_GP}
	\text{Pr}\left(||f(x)-\mu_\text{post}(x)||_2\leq \beta\sqrt{k_\text{post}(x,x)},\forall x\in\mathbb{X} \right)\geq 1-\delta,
\end{equation}
with confidence $\delta\in(0,1)$, posterior mean $\mu_\text{post}$ \eqref{post_GP_mean}, posterior variance $k_\text{post}$ \eqref{post_GP_cov}, and where $\text{Pr}({E})$ denotes the probability of an event $E$. For the scalar $\beta$, \cite{Fiedler2} (Theorem 1) proposes $\beta=M+R\sqrt{\text{log}(\text{det}(\sigma^2 I+K_X))-2\text{log}(\delta)}$ and \cite{GP_bound} (Theorem 2) $\beta=M+R\sqrt{2(\gamma+1+\text{ln}(1/\delta))}$. However, the latter requires the rather difficult calculation of the maximum information gain $\gamma$, and thus often calls for heuristic upper bounds. The error bound in \eqref{appro_GP} yields the following stochastic set membership
\begin{equation}\label{Set_mem_GP}	
\begin{aligned}
	f(x)\in\Sigma_{\text{GP}}=\{\mu(x)+R(x): \mu(x)=\mu_\text{post}(x) \text{ and }\\ R(x)^2\leq\beta^2k_\text{post}(x,x) \},
\end{aligned}
\end{equation}
with probability $1-\delta$. 

Due to the nonlinear kernel function $k$, the set membership $\Sigma_{\text{GP}}$ includes nonlinear functions as in $\Sigma_{\text{ker}}$. For that reason, \cite{Umlauft} proposes a feedback linearization to stabilize the input-affine system $\dot{x}=f(x)+G(x)u$ with a GP model for $f(x)$. The controller design requires a perfectly known and invertible input matrix $G(x)$. Moreover, a control Lyapunov function musts be calculated to establish stability guarantees, which requires the computationally complex solution of a dynamic program. Similarly, \cite{GP_Backstepping} calls for perfect knowledge of $G(x)$ and specific structure of the nonlinear dynamics $f(x)$ as common for backstepping control \citep{Khalil} (Section 13.2). Following a robust backstepping procedure, stability of the closed loop can be guaranteed from the GP inference on $f(x)$. Furthermore, \cite{GP_Lyapunov} proposes to directly study the time-derivative of a Lyapunov function for analysing the region of attraction of a stable closed loop. Indeed, the time-derivative of a Lyapunov function for a GP model of the dynamics is again a GP. Due to the impossible evaluation of the confidence intervals of this GP, a discretization of the state space and an over-\TM{estimation} by Lipschitz continuity are required, which prevent an extension to a controller synthesis. Lastly, \cite{Anne_GP}
verifies a bound on the $\mathcal{L}_2$-gain and passivity properties via optimizing over a confidence region inferred from a GP.

Since all these approaches require nonlinear optimization or a specific system dynamics, we present three approaches to tackle the nonlinearity of the kernel such that a system analysis and controller design by SDPs are possible.

\subsection{Learning linear sectors from GPs}\label{Sec_GP_Fiedler}

{\cite{Fiedler}} suggests to linearly bound the Frequentist approximation error \eqref{appro_GP}. Thereby, the nonlinear dynamics is represented by a linear system with linearly bounded uncertainty as common for system analysis and controller design by linear robust control techniques. 

For a nonlinear unknown part $\phi:[a,b]\rightarrow\mathbb{R}$ of a dynamics, \cite{Fiedler} considers the sector 
\begin{equation}\label{Sec_bound}
	\kappa_1 x^2\leq x\phi(x)\leq \kappa_2 x^2,
\end{equation} 
as common, for instance, for Lure's problem. According to the set membership $\Sigma_{\text{GP}}$ for $\phi(x)$ and \cite{Fiedler} (Lemma 2), $\phi(x)$ belongs with probability at least $1-\delta$ to the sector \eqref{Sec_bound} with 
\begin{equation}\label{Com_Sec_bound}
	\kappa_1=\min_{x\in [a,b]\backslash 0}\frac{\xi}{x^2},\ \kappa_2=\max_{x\in [a,b]\backslash 0}\frac{\xi}{x^2},
\end{equation}
$\xi=\min\left\{x(\mu_\text{post}(x)+\beta\sqrt{k_\text{post}(x,x)}), x(\mu_\text{post}(x)-\beta\sqrt{k_\text{post}(x,x)})\right\}$. Intuitively, Lemma 2 of \cite{Fiedler} determines two linear functions $\kappa_1x$ and $\kappa_2x$ that under and over approximate the nonlinear boundaries $\mu_\text{post}(x)-\beta\sqrt{k_\text{post}(x,x)}$ and $\mu_\text{post}(x)+\beta\sqrt{k_\text{post}(x,x)}$, respectively. Thereby, $\phi(x)$ is contained within the set $\{\kappa x: \kappa_1\leq\kappa\leq\kappa_2\}.$ Note that this result requires $\phi$ to satisfy Assumption~\ref{AssKernel}, i.e., $\phi$ is an element of the RKHS of the kernel of the GP and an upper bound of the norm of $\phi$ associated to the RKHS is known. Moreover, the computation of $\kappa_1,\kappa_2$ includes a nonlinear optimization problem, which might be complex for multivariate $\phi:\mathbb{R}^n\rightarrow\mathbb{R}$.

We refer to \cite{Fiedler} for the application of the linear sector \eqref{Sec_bound} for a linear state-feedback design with quadratic performance using LFRs and SDPs. Note that also other control problems can be solved based on this LFR description.

The presented concept is comparable with the polynomial representation from Section~\ref{Sec_SRNL} because both determine a suitable sector around the nonlinearity. However, \cite{MartinTP2} directly determines a polynomial sector of the dynamics instead of the inference of a learning method. Therefore, depending on the nonlinearity of the learning method, the linear sector \eqref{Sec_bound} might be more conservative than actually necessary for the nonlinear dynamics. For instance, the shape of $\mu_\text{post}(x)\pm\beta\sqrt{k_\text{post}(x,x)}$ strongly depends on the chosen kernel function and its hyperparameters. Moreover, whereas \cite{MartinTP1} computes the sector by an SDP, \eqref{Com_Sec_bound} requires to solve a nonlinear optimization problem.

\subsection{Data-driven control by linearized kernels}\label{Sec_GP_Linearization}

\cite{GP_RobCon} linearizes the posterior mean \eqref{post_GP_mean} and covariance \eqref{post_GP_cov} around an equilibrium for a linear robust controller design. Therefore, \cite{GP_RobCon} considers a Bayesian rather than a Frequentist treatment as in \cite{Fiedler}. Furthermore, we analyze the stabilization of an unknown nonlinear system by nonlinearity cancellation as proposed by \cite{Persis_kernel}. 

Similar to Section~\ref{Sec_GP_Fiedler}, {\cite{GP_RobCon}} presents a robust controller design by an LFR, but of the linearized nonlinear dynamics of $x(t+1)=f(x(t),u(t))$, deduced from a GP. Since the derivative of a GP is again a GP, \cite{GP_RobCon} infers the Jacobian matrix of $f(x,u)=\begin{bmatrix}f_1(x,u) & \cdots & f_{n_x}(x,u)\end{bmatrix}^T$ around an equilibrium point $\xi_\text{e}=\begin{bmatrix}x^T_\text{e}&u^T_\text{e}\end{bmatrix}^T$ by
\begin{equation}\label{GP_Jac}
	\frac{\partial f_i}{\partial \xi}\Bigg|_{\xi_\text{e}}\sim\mathcal{N}\left(\mu'_i(\xi_\text{e}), \Xi_i'(\xi_\text{e})\right).
\end{equation}
The posterior mean function $\mu'_i(\xi_\text{e})$ and the covariance matrix $\Xi_i'(\xi_\text{e})$ can be found in \cite{GP_RobCon} (equation (10) and (11)). Thus, \cite{GP_RobCon} can conclude on a probabilistic set membership for the Jacobian matrix $\frac{\partial f}{\partial \xi}\Big|_{\xi_\text{e}}$.
%
We refer to equations (17)-(21) in \cite{GP_RobCon} for an LFR for the linearized dynamics and to Theorem 1 in \cite{GP_RobCon} for a robust linear state-feedback synthesis with $\mathcal{H}_2$ performance by solving an SDP with LMI constraints. Moreover, \cite{GP_RobCon} suggests possible extensions by learning the operating point $\xi_\text{e}$ from the learned GP of $\Phi$, updating the GP by additional data to improve the control performance, not full state measurements, and tracking control.    

We emphasize that \cite{GP_RobCon} elaborates a Bayesian inference on the system dynamics, and therefore the data-generating $f(x,u)$ is a sample from a GP. Since the controller from \cite{GP_RobCon} (Theorem 1) robustly asymptotically stabilizes the set membership of the Jacobian matrix, the controller asymptotically stabilizes the nonlinear system if its Jacobian matrix sampled from \eqref{GP_Jac} is contained within this set. But in fact, no guarantees regarding the region of attraction and the performance of the closed loop can be deduced because the high order nonlinearities are neglected for the synthesis. Contrary, \cite{Fiedler} incorporates the nonlinearity within the synthesis, and thereby can, e.g., guarantee closed-loop performance and a region of attraction. For that reason, the procedure of \cite{GP_RobCon} is not suitable for determining dissipativity properties. Furthermore, instead of deriving an inference on the Jacobian linearization by first learning a GP, one could also directly receive the Jacobian from data following the polynomial approximation in Section~\ref{Sec_SRNL} for a TP with $\omega=\xi_\text{e}$ and $k=1$. \\      

{\cite{Persis_kernel}} presents a second approach to stabilize a nonlinear system by mainly focusing on the linear part of a kernel ridge regression. More specifically, \cite{Persis_kernel} considers the nonlinear dynamics
\begin{equation*}
	x(t+1)=f(x(t))+Bu(t),
\end{equation*}
with unknown nonlinear drift $f(x)$ and unknown input matrix $B$. To infer on the drift, data $\{{y}_i,{x}_i\}_{i=1}^{S}$ with $y_i=f(x_i)$ of the system without exciting input $u=0$ have to be measured. Applying \eqref{Set_mem_kernel} for each row implies under Assumption~\ref{AssKernel} the uncertain system representation
\begin{equation*}
	x(t+1)=AK(x(t))+Bu(t)+R(x(t)),
\end{equation*}
with data-dependent matrix $A=\begin{bmatrix}y_1 & \cdots & y_S\end{bmatrix}(\lambda I+K_X)^{-1}$ and $||R(x)||_2^2\leq n_xM^2(k(x,x)-K(x)^T\hat{K}_X^{-1}K(x))$. 

To conclude on the unknown input matrix $B$, additional data $\{\tilde{y}_i,\tilde{x}_i,\tilde{u}_i\}_{i=1}^{\tilde{S}}$ with $\tilde{y}_i=f(\tilde{x}_i)+B\tilde{u}_i$ is measured. Alternatively to the cumulative uncertainty bound \citep{Persis_kernel} (equation (31) and (32)), we suggest to compute a tighter ellipsoidal outer approximation $\Sigma_B$ for $||\tilde{y}_i-AK(\tilde{x}_i)-B\tilde{u}_i||^2_2\leq n_xM^2(k(\tilde{x}_i,\tilde{x}_i)-K(\tilde{x}_i)^T\hat{K}_X^{-1}K(\tilde{x}_i)),i=1,\dots,\tilde{S}$, as in \cite{MartinTP2}. 

To stabilize the nonlinear closed loop, \cite{Persis_kernel} separates the known kernel approximation $AK(x(t))$ into its linear and nonlinear components and uses the control structure $u=\bar{F}x+\tilde{F}\tilde{K}(x)$.
Then, \cite{Persis_kernel} follows the nonlinearity cancellation approach of \cite{NonlinearityCancell}. Hence, the linear feedback $\bar{F}x$ aims to stabilize the linear closed-loop dynamics, 
whereas the nonlinear feedback $\tilde{F}\tilde{K}(x)$ tries to minimize the influence of the nonlinearity of the closed loop. 
Both can be formulated as an SDP \citep{Persis_kernel}.

Following the discussion of \cite{GP_RobCon}, the approach of \cite{Persis_kernel} can only guarantee asymptotic stability if the choice of kernels satisfies $\lim_{||x||_2\rightarrow0}\sqrt{k(x,x)-K(x)^T\hat{K}_X^{-1}K(x)}/||x||_2 =0$. Since $R(x)$ is neglected for stabilization, weaker closed-loop guarantees are achieved compared to \cite{Fiedler}. Moreover, since the system representation and the closed-loop representation include the nonlinear term $R(x)$, the verification of dissipativity or a region of attraction would require nonlinear optimization \citep{Persis_kernel} (Theorem 2).

\subsection{Polynomial kernel for data-driven control}\label{Sec_GP_poly}

In contrast to Section~\ref{Sec_GP_Fiedler} and Section~\ref{Sec_GP_Linearization}, {\cite{Poly_Kernel}} does not reduce the nonlinear kernel to a linear representation but to a polynomial. To this end, the nonlinearity of kernels is handled by polynomial kernels together with an uniform bound for the approximation error of the nonlinear dynamics by polynomials. Due to the polynomial system representation, SOS optimization can be applied for verifying system properties and controller synthesis.

As in Section~\ref{Sec_Poly_rational}, the application of SOS relaxation motivates the investigation of continuous-time systems
\begin{equation*}
	\dot{x}(t)=f(x(t))+G(x(t))u(t)+\Phi(x(t)),
\end{equation*}  
where $f:\mathbb{X}\rightarrow\mathbb{R}^{n_x}$ and $G:\mathbb{X}\rightarrow\mathbb{R}^{n_x\times n_u}$ correspond to polynomial prior knowledge on the dynamics and $\Phi:\mathbb{X}\rightarrow\mathbb{R}^{n_x}$ to an unknown and potentially nonlinear term. Furthermore, let $f(0)=\Phi(0)=0$. To infer the unknown nonlinear function $\Phi$, let data $\{y_i,x_i,u_i\}_{i=1}^S$ with $y_i=f(x_i)+G(x_i)u_i+\Phi(x_i)+d_i$ and uniformly bounded noise $||d_i||_\infty\leq\sigma$ be given. Then \cite{Poly_Kernel} suggests to approximate each element of $\Phi=\begin{bmatrix}\Phi_1(x) & \cdots & \Phi_{n_x}(x)\end{bmatrix}^T$ by polynomial kernels 
\begin{equation*}
	k(x,x')=\sum_{i=1}^{\ell}\alpha_i^2(x^Tx')^i,
\end{equation*}
with scaling factors $\alpha_i\in\mathbb{R}$. Since the corresponding RKHS $\mathcal{H}_\text{poly}$ is the set of all polynomials of degree less than or equal to $\ell$ and zero at zero, $\Phi\notin\mathcal{H}_\text{poly}$. Thus, Assumption~\ref{AssKernel} is violated. To this end, \cite{Poly_Kernel} supposes the following assumption.

\begin{assu}[\cite{Poly_Kernel} (Assumption 3)]\label{Ass_Poly_Kern}
For a known $\epsilon>0$, there exists a polynomial vector $q(x)$ of degree less or equal to $\ell$ such that $||\Phi(x)-q(x)||_\infty\leq\epsilon$ for all $x\in\mathbb{X}$. 
\end{assu}
By Assumption~\ref{Ass_Poly_Kern}, the data satisfy
\begin{equation}\label{Data_poly_Kern}
	y_i=f(x_i)+G(x_i)u_i+q(x_i)+\tilde{d}_i,
\end{equation} 
with $||\tilde{d}_i||_\infty\leq \sigma+\epsilon$. Since each element of $q(x)=\begin{bmatrix}q_1(x) & \cdots & q_{n_x}(x)\end{bmatrix}^T$ is an element of $\mathcal{H}_\text{poly}$, we can derive for each $q_i(x)$ the set membership \eqref{Set_mem_GP} from data \eqref{Data_poly_Kern} and known $||q_i||_{\mathcal{H}_\text{poly}}$ according to Assumption~\ref{AssKernel}. This in turn results in a set membership  
\begin{equation}\label{Set_mem_poly_kern}
	\{\mu(x)+R(x): \mu(x)=\mu_\text{post}(x) \text{ and }\\ R(x)^2\leq\beta^2k_\text{post}(x,x)+\epsilon^2 \}.
\end{equation}
for each $\Phi_i(x)$.

\cite{Poly_Kernel} shows that the posterior mean $\mu_\text{post}(x)$ and variance $k_\text{post}(x,x)$ are polynomial for a polynomial kernel. Hence, the set membership \eqref{Set_mem_poly_kern} is characterized by polynomials analogously to the set membership from TPs \eqref{Set_mem_TP}. Thus, the system analysis and controller synthesis from Section~\ref{Sec_SRNL} by solving an SOS optimization problem is also possible for the set membership from polynomial kernels \eqref{Set_mem_poly_kern}. Indeed, the uncertain TP in \eqref{Set_mem_TP} reduces to the known $\mu_\text{post}(x)$ in \eqref{Set_mem_poly_kern}, while the square of the remainder in \eqref{Set_mem_TP} and \eqref{Set_mem_poly_kern} is upper bounded by a polynomial. Furthermore, Section~\ref{Sec_SRNL} considers a row-wise inference on the dynamics by TPs compatible to the kernel inference for each $\Phi_i$ by \eqref{Set_mem_poly_kern}. Nevertheless, \cite{Poly_Kernel} shows an alternating optimization over the state feedback and Lyapunov function. 


Originally, \cite{Poly_Kernel} (Proposition 3 and Section 4) does not explicitly account for the uncertainty $\epsilon$ in \eqref{Set_mem_poly_kern} and in the subsequent controller synthesis. If we include this uniform bounded uncertainty, then \eqref{Set_mem_poly_kern} is non-zero for $x=0$. Thus, an inference on stability of the origin is prevented. To circumvent this issue, we suggest to derive a polynomial bound $\epsilon(x)\geq0$ with $\epsilon(0)=0$. For instance, a suitable $\epsilon(x)$ can be calculated from the remainder formula for TPs \eqref{PolyRemainderBound} or Hermite polynomials \citep{Hermite} under Assumption~\ref{AssBoundDeri}.

Following the idea of computing $\epsilon(x)$ from the remainder formula for TPs, then $q_i(x)$ corresponds to the TP of $\Phi_i(x)$ in Assumption~\ref{Ass_Poly_Kern}. In this case, the set membership from TPs \eqref{Set_mem_TP} and from polynomial kernels \eqref{Set_mem_poly_kern} include the bound for the remainder \eqref{PolyRemainderBound}. Further, \eqref{Set_mem_TP} computes the smallest ellipse containing the set of polynomials consistent with \eqref{Data_poly_Kern}. The centre of the ellipse can be interpret as an LSE estimation. On the other hand, the posterior mean $\mu_\text{post}(x)$ corresponds to the LSE estimation for the TP $q_i$ from data \eqref{Data_poly_Kern}. 




%% file: 05_LPV.tex
\section{Data-driven control by LPV embedding}\label{Sec_LPV}

An LPV system is a linear system with system matrices that depend on a time-varying independent variable. This variable is called scheduling variable and can describe nonlinearities, time-varying system parameters, or exogenous effects. In contrast to linear time-varying systems, the scheduling variable is unknown a priori but measurable. Hence, gain scheduling control is possible where the control law depends on the state and the scheduling variable. While LPV systems establish an interesting extension of LTI systems in itself, \cite{LPV_Nonlinear} shows their potential of capturing the behavior of nonlinear systems.

For that reason, we review the data-driven treatment of a nonlinear system by combining data-driven control of LPV systems and the embedding of the nonlinear system into an LPV system. To this end, we first provide an overview on the extension of data-driven results for LTI systems to LPV systems. Afterwards, the embedding of nonlinear systems into LPV systems is discussed from a data-based perspective.

\subsection{Data-driven control for LPV systems}\label{Sec_DDLPV}

This section introduces three data-driven representations for unknown LPV systems: (i) an extension of Willems' fundamental lemma, (ii) an extension of \cite{PersisLinear}, and (iii) a set-membership approach. All three data-driven system representations ensure rigorous guarantees and system analysis and control based on SDPs. Notice that specific data-driven control schemes for LPV systems has been already investigated previously, for instance, by virtual reference feedback tuning \citep{LPV_VRFT}. See also the survey by \cite{DDLPV_survey}. \\

(i) In the behavioral framework, the dynamics of a system is characterized by its behavior that spans all possible input-output trajectories that can be observed from the system. Since the framework is trajectory-based, it is appealing for representations of dynamical systems from measured trajectories. One particular result is the fundamental lemma by \cite{Willems} that characterizes any trajectory of an LTI system based on a single measured trajectory. Therefore, this non-parametric parametrization can directly be exploited for data-driven simulation and output matching control \citep{Markovsky}, predictive control \citep{DeePC,DDMPC,DDMPC2}, system level synthesis \citep{SLS}, and the verification of dissipativity \citep{SA_Maupong, OneShot} and more general input-output properties \citep{IQC_Anne} over a data-dependent time horizon. Moreover, various extensions of the fundamental lemma, among others, linear time-varying systems \citep{Fundamental_LTV} and Wiener or Hammerstein systems \citep{Fundamental_Berberich} are investigated.  


Based on a behavioral formulation, 
{\cite{LPV_FundamentalLemma}} (Theorem 2) presents a fundamental lemma for general LPV systems. For the sake of simplicity, we only present the fundamental lemma for LPV systems with IO-representation 
\begin{equation}\label{IO_LPV}
	y(t)=\sum_{i=1}^{n_a}a_i(p(t-i))y(k-i)+\sum_{i=1}^{n_b}b_i(p(t-i))u(t-i),
\end{equation} 
with output $y(t)\in\mathbb{R}^{n_y}$, input $u(t)\in\mathbb{R}^{n_u}$, scheduling parameter $p(t)=\begin{bmatrix}p_1(t) & \cdots & p_{n_p}(t)\end{bmatrix}^T\in\mathbb{P}$
, and $n_a\geq1$, $n_b\geq0$. The functions $a_i$ and $b_i$ depend affine on the time-shifted scheduling parameter, i.e.,
\begin{align*}
	a_i(p(t-i))&=\sum_{j=0}^{n_p}a_{i,j}p_j(t-i),\\
	b_i(p(k-i))&=\sum_{j=0}^{n_p}b_{i,j}p_j(t-i),
\end{align*}
with unknown coefficients $a_{i,j}\in\mathbb{R}^{n_y\times n_y}$ and $b_{i,j}\in\mathbb{R}^{n_y\times n_u}$. Now let a trajectory $(y,p,u)$ from \eqref{IO_LPV} of length $N$ be given, i.e., $y=\begin{bmatrix}y(1) & \cdots & y(N)\end{bmatrix}$, $p=\begin{bmatrix}p(1) & \cdots & p(N)\end{bmatrix}$, and $u=\begin{bmatrix}u(1) & \cdots & u(N)\end{bmatrix}$ satisfying \eqref{IO_LPV}. Furthermore, let $(u,p)$ be PE of a sufficiently large degree \citep{LPV_FundamentalLemma}. Then the fundamental lemma for LTI systems implies that for any trajectory $(\bar{y},\bar{p},\bar{u})$ from \eqref{IO_LPV} of length $L$, there exists a $g\in\mathbb{R}^{N-L+1}$ such that
\begin{equation}\label{FL_LPV}
\begin{bmatrix}
	\mathcal{H}_L(u)\\ \mathcal{H}_L(p\otimes u)-\bar{P}_{n_u}\mathcal{H}_L(u)\\ \mathcal{H}_L(y)\\ \mathcal{H}_L(p\otimes y)-\bar{P}_{n_y}\mathcal{H}_L(y)
\end{bmatrix}g=\begin{bmatrix}\text{vec}(\bar{u})\\0\\ \text{vec}(\bar{y})\\0\end{bmatrix},
\end{equation}
where $\otimes$ is the Kronecker product of two matrices, $\bar{P}_{n}$ is a block-diagonal matrix with diagonal blocks $\bar{p}(t)\otimes I_n$, $\text{vec}(\bar{u})=\begin{bmatrix}\bar{u}(1)^T & \cdots & \bar{u}(L)^T\end{bmatrix}^T$, and $\mathcal{H}_L$ the Hankel matrix of depth $L$ \citep{LPV_FundamentalLemma}. 

By the simple algebraic relation \eqref{FL_LPV}, a single trajectory $(y,p,u)$ is sufficient to generate any trajectory $(\bar{y},\bar{p},\bar{u})$ of an LPV system. Thus, it can be leveraged to predict the behavior of an LPV system for a data-driven predictive control scheme \citep{LPV_PC}. Following the ideas for LTI systems \citep{OneShot}, \cite{LPV_dissi} determines finite-horizon dissipativity properties from one recorded noisefree trajectory. 
\\

(ii) While the fundamental lemma \eqref{FL_LPV} allows for an easy treatment of input-output data under rather mild controllability assumptions on the system, it only allows for system properties to be inferred over finite horizon. For arbitrary time horizons, {\cite{LPV_direct}} extends \cite{PersisLinear} to a representation of open- and closed-loop LPV systems from input-scheduling-parameter-state data. For that purpose, let the system matrices $A:\mathbb{P}\rightarrow\mathbb{R}^{n_x\times n_x}$ and $B:\mathbb{P}\rightarrow\mathbb{R}^{n_x\times n_u}$ be affine w.r.t. the scheduling parameter $p$
, i.e.,
\begin{equation}\label{affine_LPV}
\begin{aligned}
	A(p)=A_0+\sum_{i=1}^{n_p}p_iA_i,\ 
	B(p)=B_0+\sum_{i=1}^{n_p}p_iB_i,
	\end{aligned}
\end{equation}
with unknown matrices $A_i,B_i,i=1,\dots,n_p$. Then \cite{LPV_direct} studies an LPV system in SS-representation 
\begin{equation}\label{LPV_SS}
\begin{aligned}
	x(t+1)&=A(p(t))x(t)+B(p(t))u(t)\\
	&=\mathcal{A}\begin{bmatrix}x(t)\\ p(t)\otimes x(t)	\end{bmatrix}+\mathcal{B}\begin{bmatrix}u(t)\\ p(t)\otimes u(t)	\end{bmatrix},
\end{aligned}
\end{equation}
with state $x(t)\in\mathbb{R}^{n_x}$, input $u(t)\in\mathbb{R}^{n_u}$, and unknown matrices $\mathcal{A}=\begin{bmatrix}A_0 & A_1 & \cdots & A_{n_p}\end{bmatrix}$ and $\mathcal{B}=\begin{bmatrix}B_0 & B_1 & \cdots & B_{n_p}\end{bmatrix}$.

To find $\mathcal{A}$ and $\mathcal{B}$, let the noise-free data $\{x_i,p_i,u_i\}_{i=1}^S$ from \eqref{LPV_SS} be available and let the data-depended matrix
\begin{equation*}
	\mathcal{G}=\begin{bmatrix}
		x_1 & \cdots & x_{S-1}\\
		p_1\otimes x_1 & \cdots & p_{S-1}\otimes x_{S-1}\\
		u_1 & \cdots & u_{S-1}\\
		p_1\otimes u_1 & \cdots & p_{S-1}\otimes u_{S-1}
	\end{bmatrix}
\end{equation*}
have full row rank. The latter condition corresponds to a PE condition that generalizes the corresponding condition for LTI systems \citep{PersisLinear}. Since 
\begin{equation*}\label{LPV_Dyn_data}
	X^+=\begin{bmatrix}x_2 & \cdots & x_S\end{bmatrix}=\begin{bmatrix}\mathcal{A} & \mathcal{B}\end{bmatrix}\mathcal{G}
\end{equation*}
and $\mathcal{G}$ has full row rank, we obtain the unknown matrices $\mathcal{A}$ and $\mathcal{B}$ by $\begin{bmatrix}\mathcal{A} & \mathcal{B}\end{bmatrix}=X^+\mathcal{G}^\dagger$ and the following data-based representation of \eqref{LPV_SS} 
\begin{equation*}
	x(t+1)=X^+\mathcal{G}^\dagger\begin{bmatrix}x(t)\\ p(t)\otimes x(t)	\\u(t)\\ p(t)\otimes u(t)	\end{bmatrix}.
\end{equation*}
By this data-based representation and the model-based conditions for stability and quadratic performance for LPV systems \citep{LPV_direct} (Lemma 1-4), we can analyze an LPV system directly from data.

For the data-driven controller synthesis for \eqref{LPV_SS}, \cite{LPV_direct} extends the data-based closed-loop representation of LTI systems \citep{PersisLinear} (Theorem 2), as shown in Section~\ref{Sec_Direct}, for a state feedback $u(t)=K(p(t))x(t)$. For a polytopic and convex $\mathbb{P}$, the design of a stabilizing feedback boils down to a finite number of LMI constraints by the full-block S-procedure \citep{SProcedure}. \cite{LPV_direct} (Theorem 4-6) presents further statements for a data-driven controller synthesis with performance guarantees. Thereby, it is possible to design a controller for LPV systems with stability and performance guarantees from data by solving an SDP with a finite number of LMIs constraints.\\

(iii) \cite{LPV_set_membership} considers a set-membership approach for LPV systems, and thus ensures rigorous guarantees despite noise-corrupted data. Moreover, the controller synthesis does not scale with the number of data as in \cite{LPV_direct}. To this end, noisy samples $\{x_i,p_i,u_i\}_{i=1}^S$ from the LPV systems in SS-representation \eqref{LPV_SS} with $B(p)=B$ are drawn, i.e.,
\begin{equation*}
	x_{i+1} =\mathcal{A}\begin{bmatrix}x_i\\ p_i\otimes x_i	\end{bmatrix}+Bu_i+d_i.
\end{equation*}
The actual noise realizations $d_i,i=1,\dots,S,$ are unknown. Collecting the data into matrices $X=\begin{bmatrix}x_1 & \cdots & x_{S-1}\end{bmatrix},U=\begin{bmatrix}u_1 & \cdots & u_{S-1}\end{bmatrix}, X^+=\begin{bmatrix}x_2 & \cdots & x_{S}\end{bmatrix},P_X=\begin{bmatrix}p_1\otimes x_1 & \cdots & p_{S-1}\otimes x_{S-1}\end{bmatrix},$ and $D=\begin{bmatrix}d_1 & \cdots & d_{S-1}\end{bmatrix}$ amounts to the relation 
\begin{equation}\label{Data_LPV}
	X^+=\begin{bmatrix}\mathcal{A}& B\end{bmatrix}\begin{bmatrix}X\\ P_X\\ U	\end{bmatrix}+D.
\end{equation}
For simplicity and analogously to the LTI result \citep{vanWaarde}, let the noise be cumulatively bounded
\begin{equation*}\label{Energy_noise_bound}
	\begin{bmatrix}I\\ D^T\end{bmatrix}^T\begin{bmatrix}\Delta_1 & \Delta_2^T\\ \Delta_2 & \Delta_3\end{bmatrix}\begin{bmatrix}I\\ D^T\end{bmatrix}\succeq0,
\end{equation*}   
with $\Delta_3\preceq0$. See \cite{MartinIQC} for a comparison of this cumulative and pointwise-in-time noise bound. Substituting the unknown noise matrix $D$ from \eqref{Data_LPV} immediately yields the set membership for the unknown coefficient matrices $\mathcal{A}$ and $B$
\begin{equation}\label{Set_mem_LPV}
\begin{aligned}
&	\Sigma_{\mathcal{A}, B}{=}\left\{\begin{bmatrix}\tilde{\mathcal{A}}& \tilde{B}\end{bmatrix}:\right.\\
	&\hspace{0.5cm}\left.\star^T\begin{bmatrix}\Delta_1 & \Delta_2^T\\ \Delta_2 & \Delta_3\end{bmatrix}\cdot\begin{bmatrix}I & 0\\ X^{+^T} & -\begin{bmatrix}X\\ P_X\\ U	\end{bmatrix}^T 
	\end{bmatrix}\begin{bmatrix}I\\\begin{bmatrix}\tilde{\mathcal{A}}& \tilde{B}\end{bmatrix}^T
	\end{bmatrix}\succeq0\right\}.
\end{aligned}
\end{equation} 
By ``$\star$'', the matrices on the right-hand side of ``$\cdot$'' are abbreviated.

Since $\Sigma_{\mathcal{A}, B}$ characterizes the set of all LPV systems consistent with the noisy data, the ground-truth LPV system \eqref{LPV_SS} is contained within the set. Therefore, finding a control policy, that stabilizes all LPV systems within $\Sigma_{\mathcal{A}, B}$, also stabilizes the ground-truth system. In contrast to the previous LPV approaches, here the data is corrupted by noise, and thereby the true system can not be exactly identified from data. This uncertainty adds conservatism to the controller design but can be handled by robust control techniques.


By Lemma 1 from \cite{LPV_set_membership}, the controller synthesis for LPV systems with convex polytopic scheduling space $\mathbb{P}$ boils down to the robust stabilization of the LTI systems with system matrices contained in the set membership \eqref{Set_mem_LPV} and scheduling parameters at the vertices $\mathbb{P}$. This can be solved by a non-conservative S-Lemma \citep{vanWaarde}. Thereby, also a controller design with performance guarantees is conceivable. In contrast to \cite{LPV_direct}, the set-membership approach has the advantage that the number of optimization variables does not increase with the number of samples.

Following a similar procedure, \cite{LPV_SetMem2} presents a set-membership approach for LPV systems using pointwise-in-time bounded noise. Moreover, the synthesis of a gain-scheduling feedback controller that ensures a robust invariant set w.r.t. bounded disturbance is investigated.


\subsection{Data-driven LPV embedding of nonlinear systems}\label{Sec_LPVEmb}

While the data-driven LPV methods from the previous section are also of independent interest, we clarify next their applicability for checking system properties and controller design for nonlinear systems. To this end, we consider the two conversions of a nonlinear system into an LPV system as proposed by \cite{LPV_nonlinear} and \cite{LPV_nonlinear2} in a data-driven control context.

{\cite{LPV_nonlinear}} studies a general nonlinear discrete-time system
\begin{equation}\label{NL_sys_LPV}
	x(t+1)=f(x(t),u(t)),
\end{equation}
with continuously differentiable function $f:\mathbb{X}\times\mathbb{U}\rightarrow\mathbb{X}$ and $f(0,0)=0$. As shown in \cite{Extended_lin} (Proposition 1) and the references therein, the nonlinear system can be embedded into the LPV system \eqref{LPV_SS} with scheduling parameter $p(t)=\psi(x(t),u(t))\in\mathbb{P}$. Thus, $f(x,u)\in\{A(p)x+B(p)u:p\in\mathbb{P}\}$.


\begin{assu}[\cite{LPV_nonlinear}]\label{Assu_LPV}
Let the scheduling map $\psi(x,u)$ be known and chosen such that the matrices $A(p)$ and $B(p)$ of \eqref{LPV_SS} are affine in $p$ as in \eqref{affine_LPV}.	 
\end{assu}

Under Assumption~\ref{Assu_LPV}, input-scheduling-parameter-state data can be obtained from recorded input-state trajectories of the nonlinear system \eqref{NL_sys_LPV} as required for the three data-driven representations for LPV from Section~\ref{Sec_DDLPV}. Thus, a set membership for the nonlinear system \eqref{NL_sys_LPV} is given by the set of all LPV systems consistent with the data, i.e., 
\begin{equation}\label{Set_mem_LPV1}
\Sigma_{\text{LPV}}=\left\{\mathcal{A}\begin{bmatrix}x\\ p\otimes x	\end{bmatrix}+\mathcal{B}\begin{bmatrix}u\\ p\otimes u	\end{bmatrix}:p\in\mathbb{P}, \begin{bmatrix}\mathcal{A} & \mathcal{B}\end{bmatrix}=X^+\mathcal{G}^\dagger\right\}
\end{equation} 
or 
\begin{equation}\label{Set_mem_LPV2}
\tilde{\Sigma}_{\text{LPV}}=\left\{\tilde{\mathcal{A}}\begin{bmatrix}x\\ p\otimes x	\end{bmatrix}+\tilde{\mathcal{B}}\begin{bmatrix}u\\ p\otimes u	\end{bmatrix}:p\in\mathbb{P}, \begin{bmatrix}\tilde{\mathcal{A}} & \tilde{\mathcal{B}}\end{bmatrix}\in\Sigma_{\mathcal{A}, B}\right\}.
\end{equation} 
Moreover, Assumption~\ref{Assu_LPV} allows the scheduling parameter to be constructed from input-state measurements of the nonlinear system to realize the feedback law $u(t)=K(p(t))x(t)$. Hence, we can conclude that all three data-driven LPV approaches from Section~\ref{Sec_DDLPV} are conceivable for a data-driven system analysis or controller design for nonlinear systems under Assumption~\ref{Assu_LPV}. Note that \cite{LPV_nonlinear} focuses on the LPV representation (ii) from \cite{LPV_direct}.\\ 


As shown in \cite{LPV_embedding}, the previous direct LPV embedding comes with the shortcoming that stability from the LPV system \eqref{LPV_SS} implies stability for the nonlinear system \eqref{NL_sys_LPV} only at the origin but not for potential non-zero equilibria. This phenomenon emerges as the LPV stability analysis supposes that the scheduling parameter is an exogenous variable independent of state and input. Therefore, {\cite{LPV_nonlinear2}} considers an LPV embedding of the so-called velocity-form of a nonlinear system. Thereby, global stability and performance guarantees are possible.    

The velocity-form of a nonlinear system describes the dynamics of $\Delta x(t)=x(t)-x(t-1)$ and naturally exhibits a state-depended form. 
Similar to Assumption~\ref{Assu_LPV}, \cite{LPV_nonlinear2} requires additional insights into the nonlinearities of its dynamics to infer a suitable scheduling variable.
%
As for the direct LPV embedding, the data-driven control methods for LPV systems (i)-(iii) from the previous section are applicable for analysing and designing controllers for the embedded velocity-form LPV system using SDPs with a finite number of LMI constraints. 
Also note the similarities of the analysis of the velocity-form and incremental system analysis of nonlinear system \eqref{NL_sys_LPV} \citep{DifferentialPass}.\\ 

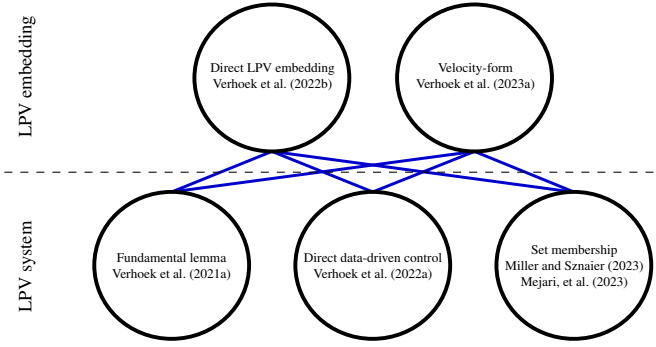
\begin{figure}
	\centering
	\begin{tikzpicture}
	\node[inner sep=0pt] (p1) at (0,0) {\includegraphics[width=.45\textwidth]{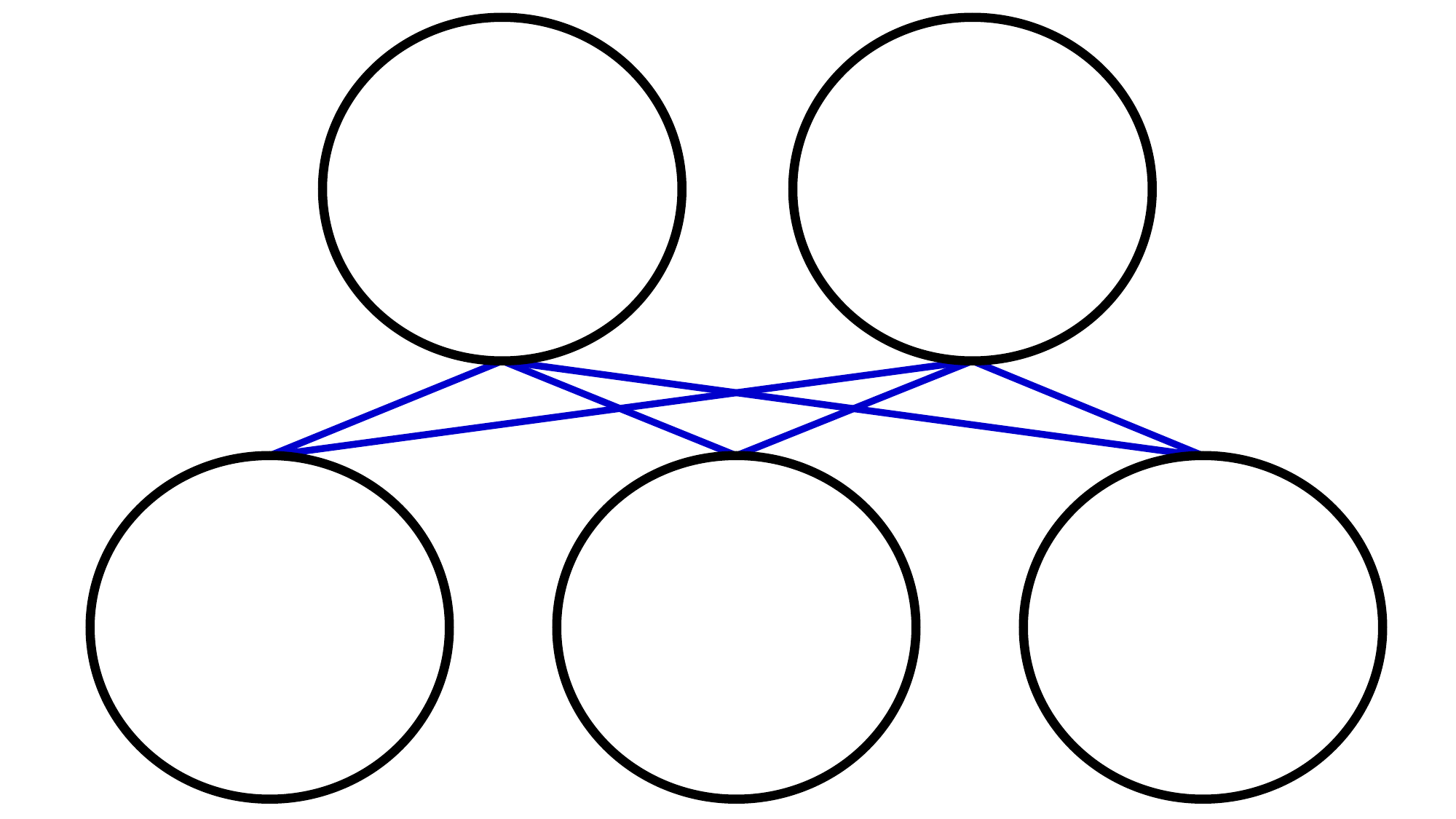}};
	
	\node[inner sep=0pt] (p1) at (-2.6,-1.25) {\tiny \begin{tabular}{c} Fundamental lemma \\ \cite{LPV_FundamentalLemma} \end{tabular} };
	\node[inner sep=0pt] (p1) at (0.05,-1.25) {\tiny \begin{tabular}{c} Direct data-driven control \\ \cite{LPV_direct} \end{tabular} };
	\node[inner sep=0pt] (p1) at (2.7,-1.25) {\tiny \begin{tabular}{c} Set membership \\ \cite{LPV_set_membership}\\ \cite{LPV_SetMem2} \end{tabular} };

	\node[inner sep=0pt,rotate=90] (p1) at (-4.5,-1.25) {\scriptsize  LPV system};
	\node[inner sep=0pt,rotate=90] (p1) at (-4.5,1.27) {\scriptsize LPV embedding};

	\node[inner sep=0pt] (p1) at (-1.27,1.27) {\tiny \begin{tabular}{c} Direct LPV embedding \\ \cite{LPV_nonlinear} \end{tabular} };
	\node[inner sep=0pt] (p1) at (1.4,1.27) {\tiny \begin{tabular}{c} Velocity-form \\ \cite{LPV_nonlinear2} \end{tabular} };
	
	\draw[dashed] (-4.8,0)--(3.8,0);
	\end{tikzpicture}
	\caption{Combining data-driven control techniques for LPV systems and LPV embedding of nonlinear systems.} 
	\label{Fig.LPV_summary}
\end{figure}

As summarized in Figure~\ref{Fig.LPV_summary}, combining an embedding of an unknown nonlinear dynamics into an LPV system and data-driven control techniques for LPV systems enables a data-based system analysis and control of nonlinear systems. 
The key idea to analyze systems and design controllers by SDPs despite nonlinear dynamics is the linearization by embedding the nonlinear system into an LPV system. Thereby, the nonlinearity is relaxed as the free scheduling variable. The LPV embedding comes at the cost of additional conservatism due to the independence of the scheduling variable regarding states and inputs. Thus, the LPV representation allows for more possible trajectories than the ground-truth nonlinear system. Nevertheless, \cite{LPV_Nonlinear} indicates the successful application of an LPV embedding for a large subset of nonlinear systems. We also refer to \cite{LTI_embedding} for an embedding of a nonlinear system into an LTI system by introducing additional inputs and for its application for data-driven simulation.    

In contrast to the previously presented approaches by polynomial approximation (Section \ref{Sec_PolyApprox}) and kernel regression (Section \ref{Sec_GP}), the LPV approach requires knowledge of a function basis containing the nonlinearity of the dynamics (Assumption~\ref{Assu_LPV}) or of its velocity-form. This insight might be reasonable for mechanical or electrical systems, however, is often not known. The latter scenario is tackled in the polynomial approximation and kernel regression approach. We assess the latter problem more difficult as the data only provide local insights into the nonlinear dynamics. Hence, the number of required data samples, e.g., in \cite{LPV_nonlinear}, can be significantly smaller. Further advantages of the LPV approach are the design of a potentially more flexible nonlinear feedback law and it does not rely on the scalability of an SOS relaxation. 

We conclude with some possible extensions and open questions. As commented in \cite{LPV_nonlinear2} (Remark 1), Assumption~\ref{Assu_LPV} could be relaxed by using a monomial basis instead of a true function basis. This motivates the question whether an upper bound on the remainder for TPs and the insight from Assumption~\ref{AssBoundDeri} can be incorporated into the data-based LPV representations and the subsequent analysis. For a monomial basis, a second question is whether a direct SOS treatment of the polynomial approximation is beneficial compared to an LPV relaxation by embedding the polynomial dynamics into an LPV system. Furthermore, we point out that an LPV embedding under Assumption~\ref{Assu_LPV} might be non-unique. For instance, given $p_1=x_1x_2$ and $p_2=x_2^2$, then
\begin{equation*}
	\begin{bmatrix}x_1^+\\x_2^+
	\end{bmatrix}=\begin{bmatrix}x_1x_2^2\\x_1^2x_2+x_2^3\end{bmatrix}=\begin{bmatrix}\alpha p_2 & (1-\alpha)p_1\\ p_1 & p_2\end{bmatrix}\begin{bmatrix}x_1\\x_2
\end{bmatrix} \forall\alpha\in\mathbb{R}.
\end{equation*}
Hence, the problem arises that the set of matrices $A_1$ and $A_2$ from \eqref{affine_LPV} compatible with the data is always a linear subspace, and thus unbounded. At the same time, many robust control techniques require a bounded uncertainty set.

\subsection{Iterative control scheme by extended linearization}\label{Sec_SysRepExLin}

Related to the problem setup in \cite{LPV_nonlinear}, we briefly report the data-driven iterative control scheme of {\cite{Sznaier}}. There the controller design for an unknown nonlinear system $x(t+1)=f(x(t))+g(x(t))u(t)$ is solved by the extended linearization \citep{Extended_lin} and under a known function basis that spans the dynamics. 
  	
\begin{assu}[\cite{Sznaier} (Assumption 1)]\label{Assu_ExtLin}
Let basis functions $F:\mathbb{R}^{n_x}\rightarrow\mathbb{R}^{n_f}$ and $G:\mathbb{R}^{n_x}\rightarrow\mathbb{R}^{n_g\times n_u}$ be known such that there exist matrices $A\in\mathbb{R}^{n_x\times n_f},B\in\mathbb{R}^{n_x\times n_g}$ with $f(x)=A F(x)$ and $g(x)=B G(x)$.
\end{assu}

Assumption~\ref{Assu_ExtLin} calls for a dictionary of basis functions for the nonlinear system dynamics. This insight might by accessible from first principles, e.g., for mechanical or electrical systems. 

Proceeding as in \cite{LPV_set_membership} or \cite{LPV_SetMem2} yields a set membership $\Sigma_{A,B}$ for the unknown matrices $A$ and $B$ from noisy data $\{x_i,u_i\}_{i=1}^S$. Thus, \cite{Sznaier} concludes on the nonlinear data-based system representation
\begin{equation}\label{DD_ext_lin}
	f(x)+g(x)u\in\left\{\tilde{A} Z(x)x+\tilde{B}G(x)u: \begin{bmatrix}\tilde{A}& \tilde{B}\end{bmatrix}\in\Sigma_{A,B}\right\},
\end{equation}
with extended linearization $F(x)=Z(x)x$.

Instead of exploiting an LPV embedding to linearize the nonlinearity of \eqref{DD_ext_lin} as in \cite{LPV_nonlinear}, \cite{Sznaier} suggests an online scheme. There the optimization problem
\begin{equation}\label{Online_problem}
\begin{aligned}
	&\min_{u(t),u(t+1),\dots}\sum_{i=t}^{\infty}x(i)^TQx(i)+u(i)^TRu(i)\\
	&\hspace{0.5cm}\text{s.t.}\ x(i+1)=\tilde{A} Z(x(t))x(i)+\tilde{B}G(x(t))u(i)
\end{aligned}
\end{equation} 
is solved in each time instant $t$ and the first part of the optimal control policy $u(t)$ is applied to the system. Hence, the key idea to circumvent the nonlinearity of \eqref{DD_ext_lin} is to freeze $Z(x)$ and $G(x)$ at time $t$ and to treat the nonlinear dynamics as an LTI system. Hence, the optimal control problem \eqref{Online_problem} can be solved for all $\begin{bmatrix}\tilde{A}&\tilde{B}\end{bmatrix}\in\Sigma_{A,B}$ by an SDP with LMI constraints using data-driven control for LTI systems \citep{vanWaarde}. To guarantee that the equilibrium at the origin is globally asymptotically stable, the optimal control problem is extended by a decrease of a control Lyapunov function along the closed-loop trajectory.   

On the one hand, this procedure constitutes an alternative to an LPV embedding as similar assumptions on prior insights and data are required. On the other hand, the optimal control problem \eqref{Online_problem} has to be solved iteratively online and the closed loop does not necessarily satisfy a specified control performance. An open question is whether the non-unique extended linearization might lead to performance or feasibility issues although not observed in the numerical example of \cite{Sznaier}. Moreover, the feasibility of the optimal control problem is not guaranteed during runtime, which is typically ensured in MPC. We highlight the connection to the data-driven predictive control approach for nonlinear systems from \cite{BerberichNL}. There also an online updated local linear approximation of the nonlinear system is exploited but by using the fundamental lemma by \cite{Willems} rather than a set-membership representation.

%% file: 06_Koopman.tex
\section{Data-driven control by state lifting}\label{Sec_Koopman}

The nonlinearity of a dynamical system \TM{usually} prevents a system analysis and controller design by convex optimization. The literature on the Koopman operator presents a solution by lifting the states to higher, potentially infinite, dimensions. Thereby, the nonlinear dynamics is exactly described by a \mbox{(bi-)linear} system, and thus can be handled by control techniques for \mbox{(bi-)linear} systems. Due to its success in application, we review the Koopman operator with a focus on results providing rigorous guarantees.

\subsection{Koopman operator paradigm}\label{Sec_Koop}

The Koopman operator paradigm was first introduced by \cite{Koopman} and it has been increasingly used over the last decades because it provides a global \mbox{(bi-)linear} system description in contrast to a local Jacobian linearization. Its applications cover prediction \citep{Mezic}, global stability analysis \citep{Mauroy}, MPC \citep{Korda}, linear-quadratic regulation \citep{Brunton}, and robotics \citep{Bruder}. Further control-oriented applications can be found in Section 5.2 of \cite{Koopman_Survey}.

The notion of the Koopman operator is originally introduced for an autonomous system $x(t+1)=f(x(t))$ with state $x(t)\in\mathbb{X}\subseteq\mathbb{R}^{n}$. Instead of observing the time evolution of the states, the Koopman operator considers the system dynamics through the lens of scalar functions $\psi:\mathbb{X}\rightarrow\mathbb{C}$ from a function space $\mathcal{H}$. Thereby, the Koopman operator $\mathcal{K}:\mathcal{H}\rightarrow\mathcal{H}$ is given by $\mathcal{K}\psi=\psi\circ f$, with function composition $\circ$. Thus, the operator characterizes the propagation of a whole hypersurface over one time-step rather than of a single state. Due to the linearity of the function composition, $\mathcal{K}$ is linear, 
although the underlying dynamics is nonlinear. Moreover, $\mathcal{K}$ operates globally for all $x\in\mathbb{X}$ with $\mathcal{K}\psi(x)=\psi\circ f(x)=\psi(x(t+1))$. 

Besides providing a global linearization, the Koopman theory can be useful in system analysis: The particular choice of observables by the eigenfunctions $\phi$ of $\mathcal{K}$, with $\mathcal{K}\phi=\lambda\phi$ and eigenvalue $\lambda$, enables a spectral analysis of nonlinear systems. Furthermore, \cite{Contaction_Koopman} investigates the equivalence of Koopman theory and contraction theory and \cite{Koopman_Lyapunov} shows the connection of Koopman and Lyapunov theory by interpreting Lyapunov functions as special case of observable. We refer to the survey by \cite{Koopman_Survey} and the references therein for more details including data-based inferences on eigenfunctions (Section 3.4).  

The linearity of the system representation comes at the cost that the Koopman operator is infinite-dimensional. Therefore, a finite dimensional truncation is necessary for an efficient analysis and controller design. For that purpose, a dictionary of finitely many observables $\{\psi_i\}_{i=1}^{n_D}$ is chosen, which leads to the lifted finite-dimensional dynamics $z(t+1)=Az(t)$ with $z=\begin{bmatrix}\psi_1(x) & \cdots & \psi_{n_D}(x)\end{bmatrix}^T=\Psi(x)$. However, since the dictionary generally does not span an invariant space w.r.t. the Koopman operator, the lifted dynamics involves a truncation error. Since this error depends on the system dynamics and the chosen observables, its structure and size is in general unclear. The literature examines, besides specific choices of observables as time-delay coordinates \citep{TimeDelayKoopman}, learning of suitable observables from data, e.g., learning the eigenfunctions \citep{Learning_Koopman_Eigenfcn}, deep learning \citep{Koopman_NN}, and by SDPs \citep{Koopman_sznaier}. While the observables are mostly chosen to achieve a small prediction error, optimizing the dictionary together with a controller is interesting to find observables achieving good control performance.

While the Koopman operator is linear for autonomous systems, the Koopman paradigm leads for input-affine systems and state-dependent lifting functions to a bilinear system description in continuous-time and to an LPV description in discrete-time \citep{Koopman_Survey} (Corollary 5.1.1). Nevertheless, most works restrict the finite-dimensional truncation model to an LTI system \citep{Korda, Koopman_Guarantees} $z(t+1)=Az(t)+Bu(t)$ or a bilinear system \citep{Koopman_bilnear, Strasser_Koopman} $z(t+1)=Az(t)+u(t)(B_0+B_1z(t))$ (here $u(t)\in\mathbb{R}$).  

Extended dynamic mode decomposition (EDMD) constitutes a common estimation of the finite-dimensional system matrices of the truncation models from data: For a given set of state-dependent observables $\{\psi_i\}_{i=1}^{n_D}$, input-state samples $\{x_i,u_i\}_{i=1}^S$ are collected from the underlying nonlinear system $x(t+1)=f(x(t))+g(x(t))u(t)$. Then, for a bilinear system representation, the LSE problem 
\begin{equation}\label{EDMD}
	\min_{A,B_0,B_1}\Big|\Big|Z^+-\begin{bmatrix}A & B_0 & B_1\end{bmatrix}Y\Big|\Big|_\text{Fr}
\end{equation}
is solved with data-dependent matrices $Z^+=\begin{bmatrix}\Psi(x_2) & \cdots & \Psi(x_{S})\end{bmatrix}$ and 
\begin{equation*}
	Y=\begin{bmatrix}\Psi(x_1) & \cdots & \Psi(x_{S-1})\\ u_1 & \cdots & u_{S-1}\\ \Psi(x_1)u_1 & \cdots & \Psi(x_{S-1})u_{S-1}\end{bmatrix}.
\end{equation*}
As shown in \cite{Korda2}, EDMD is asymptotically consistent, i.e., converges to the Koopman operator for $n_D\rightarrow\infty$ and $S\rightarrow\infty$ under some additional assumptions. In the non-asymptotic case, \cite{Koopman_error} provides \TM{statistical} bounds on the \TM{estimation} error due to finite data. 

\subsection{Data-driven control by state lifting with guarantees}\label{Sec_Koop_ctrl}

In the following, we take a closer look at two articles 
because they incorporate an \TM{estimation} error into the controller design by a Koopman-lifted system representation. Thereby, guarantees for the closed loop of the original nonlinear system can be recovered if the \TM{estimation} error satisfies the assumed error characterization. 

\cite{Koopman_Guarantees} proposes to lift the nonlinear system $x(t+1)=f(x(t),u(t)), x\in\mathbb{X},u\in\mathbb{U}$, for a suitable choice of observables to the LTI system
\begin{equation}\label{LTI_Koopman}
\begin{aligned}
	z(t+1)&=Az(t)+Bu(t)+w(t),\\
	{x}(t)&=Cz(t)+v(t).
\end{aligned} 
\end{equation}
The matrices $A$ and $B$ are calculated from EDMD, analogously to \eqref{EDMD} with an additional regularization, and the output matrix $C$ from 
\begin{equation*}\label{EDMD2}
	\min_{C}\sum_{i=1}^{S}||C\Psi(x_i)-x_i||_2^2+\beta\,||C||^2_\text{Fr},
\end{equation*}
with regularization parameter $\beta>0$. Thereby, the transformation from $z$ to the predicted state ${x}$ is also linearized for general observables. To account for the truncation error, EDMD for finite data, and the simplified LTI representation of the infinite dimensional Koopman model, \cite{Koopman_Guarantees} extends the lifted system by bounded uncertainties $w$ and $v$. 

\begin{assu}[\cite{Koopman_Guarantees}]\label{Assu_Koopman1}
For \eqref{LTI_Koopman}, $w\in\mathcal{W}$ and $v\in\mathcal{V}$ with known bounded sets $\mathcal{W}$ and $\mathcal{V}$.
\end{assu}
Under reasonable assumptions on the observables, \cite{Koopman_Guarantees} (Proposition 2) proves that the uncertainty sets $\mathcal{W}$ and $\mathcal{V}$ are bounded. Further, since the knowledge of $\mathcal{W}$ and $\mathcal{V}$ is non-trivial, \cite{Koopman_Guarantees} suggests to validate a choice of $\mathcal{W}$ and $\mathcal{V}$ using statistical learning theory. 

Under Assumption~\ref{Assu_Koopman1}, the nonlinear dynamics is captured by
\begin{equation*}
	f(x,u)\in\{C(A\Psi(x)+Bu+w)+v:w\in\mathcal{W},v\in\mathcal{V}\},
\end{equation*}
with $A,B,$ and $C$ from EDMD. The nonlinearity of this set membership by $\Psi(x)$ can be circumvented by analysing the lifted system \eqref{LTI_Koopman} with lifted state vector $z=\Psi(x)$. 

Since the error characterization of Assumption~\ref{Assu_Koopman1} does not vanish at the origin, this uncertainty description is not suitable for verifying dissipativity or a state-feedback design. Instead, \cite{Koopman_Guarantees} follows a robust tube-based MPC approach \citep{Robust_MPC}. Thereby, \cite{Koopman_Guarantees} proves closed-loop robustness w.r.t. the \TM{estimation} error and point-wise convergence (Theorem 3). Furthermore, an MPC scheme with nonconvex optimization is circumvented by the linear Koopman prediction model.\\ 

In contrast to a linear MPC, \cite{Strasser_Koopman} designs a state feedback from a lifted state model including a finite-gain bounded \TM{estimation} error. Since a lifting of an input-affine nonlinear system does not lead to an LTI system \citep{Koopman_Survey}, \cite{Koopman_bilnear} and \cite{Strasser_Koopman} consider a discrete-time bilinear representation 
\begin{equation*}
	z(t+1)=Az(t)+u(t)(B_0+B_1z(t))
\end{equation*}
of a finite-dimensional lifted system with scalar input $u(t)\in\mathbb{R}$. Thus, a higher accuracy and better control performance can be achieved. \cite{Koopman_bilnear} does not account for the error by \TM{estimating} the Koopman operator. Contrary, \cite{Strasser_Koopman} supposes a finite-gain bound for an additive \TM{estimation} error during closed-loop operation with $u(t)=k^Tz(t)$.
\begin{assu}[\cite{Strasser_Koopman} (Assumption 2)]\label{Assu_Koopman}
An additive \TM{estimation} error $\epsilon(z)\in\mathbb{R}^{n_D}$ satisfies a finite-gain bound $||\epsilon(z)||_2\leq L||z||_2$ with known $L\geq0$.	 
\end{assu} 
\cite{Strasser_Koopman} proposes to estimate $L$ by first approximating the Lipschitz constant of the nonlinear dynamics from data. Together with the matrices $A,B_0,B_1$ from EDMD, this results in a Lipschitz constant of $\epsilon(z)$ satisfying the finite-gain bound in Assumption~\ref{Assu_Koopman}. Alternatively, a validation procedure by Hoeffding's inequality or learning a kernel approximation for $\epsilon(z)$ to obtain a linear sector by \cite{Fiedler} are conceivable. However, since the gain $L$ is required for the closed loop and the controller gain $k$ is undetermined at first, closed-loop data is not available. Therefore, the estimation of $L$ is non-trivial and might require an iteration between controller synthesis and estimation of $L$. 

Under Assumption~\ref{Assu_Koopman}, the nonlinear dynamics of the closed loop $f(x)+g(x)k^T\Psi(x)$ is contained in
\begin{equation*}\label{Set_mem_Koopman}
\begin{aligned}
	&\{\Psi^{-1}(A\Psi(x)+k^T\Psi(x)(B_0+B_1\Psi(x))+\epsilon(\Psi(x))):\\
	&\hspace{4cm}||\epsilon(\Psi(x))||_2\leq L||\Psi(x)||_2\}.	
\end{aligned}
\end{equation*}
For the lifted state $z=\Psi(x)$, the nonlinearity of the set membership boils down to a bilinear system with finite-gain bounded uncertainty. Moreover, the inverse mapping $x=\Psi^{-1}(z)$ might be a simple matrix multiplication for certain lifting, e.g., delay or monomial coordinates \citep{Strasser_Koopman}. Thereby, a stabilizing linear feedback of the lifted states can be received by LMI robust control techniques for bilinear systems that also stabilizes the nonlinear system. 

\TM{Motivated by the bilinear model deduced from lifting techniques, we summarize in the following remark three data-driven control techniques for bilinear systems.}

\begin{rmk}[Data-driven control of bilinear systems]
\TM{For bilinear systems with unknown system matrices, \cite{Bilinear2} presents an extension of Willems' fundamental lemma by interpreting the bilinearity as an independent input. 
Similar to \cite{PersisLinear}, \cite{Bilinear} presents a data-based closed-loop parametrization in order to design a state feedback by LMIs. 
Since the bilinearity is actually known and available for feedback, \cite{Strasser_Bilinear} presents a gain-scheduling controller to achieve a better control performance than \cite{Strasser_Koopman}.}
\end{rmk}


\cite{Koopman_Guarantees} and \cite{Strasser_Koopman} consider a closed-loop analysis in the lifted states without ensuring that the lifted states actually have a preimage in the original state space. By incorporating a projection of the lifted states back on the original space together with the error bounds from \cite{Koopman_error}, \cite{Koopman_MPC} presents a nonlinear predictive controller with practical asymptotic stability guarantees. For that purpose, a similar error bound as in Assumption~\ref{Assu_Koopman} is deduced from \cite{Koopman_error} but in the original state space under the assumption of a finite and invariant dictionary of observables.\\

\TM{In contrast to the Koopman operator, Carleman lifting \citep{Carleman_original} provides by monomial observables a systematic approach to bound its estimation error from finite data. For instance, \cite{Error_bound_Carleman} and \cite{Carleman} derive a guaranteed, but potentially conservative, error bound between trajectories of an autonomous nonlinear system and the data-driven inference of a truncated lifted linear system. Hence, further examinations on these error bounds are necessary for the application, e.g., for a linear MPC design \citep{Carleman_MPC}. Furthermore, \cite{Taylor_Carleman} investigates the connection of Carleman lifting and TP approximation but neglects the error regarding the vector field of the truncated Carleman lifting. Similarly, \cite{SOS_Koopman} exploits a combination of Koopman lifting with polynomial observables and SOS optimization but also neglects the error due to finite data and truncation. However, due to the polynomial dictionary, an investigation of the \TM{estimation} error might be possible.}\\

Related to the state lifting by the Koopman operator, a nonlinear system $\dot{x}=f(x)+g(x)u$ can be polynomialized if the dynamics contain only certain nonlinearities. More specifically, the time derivatives of the nonlinear functions in $f$ and $g$ have to be written as terms of the same functions. See Table 2 in \cite{Strasser} for a collection of such functions. In this case, the nonlinear dynamics can be formulated as a finite-dimensional polynomial system $\dot{z}(t)=AZ(z(t))+BH(z(t))u(t)$ by introducing a lifted state $z$ together with some conservatism \citep{Strasser} (Remark 9). In a data-driven context, \cite{Strasser} examines polynomialization of nonlinear systems with known function basis of $f$ and $g$. Hence, the data-driven system analysis and controller design boils to the problem of polynomial systems with unknown coefficient matrices $A$ and $B$ as in Section~\ref{Sec_Poly_rational}.\\

The Koopman operator constitutes a powerful tool to analyze general nonlinear system by methods from linear or bilinear control theory. While the operator is infinite-dimensional, the literature asks for finite-dimensional \TM{estimations} retrieved from a finite set of data. Recent works take the error by truncation and simplified lifted model into account and handle it by robust control techniques. Thus, a system analysis and a controller design with rigorous guarantees \TM{and using} SDPs are possible. But in fact, the characterization of this error is so far rather vague, and therefore potentially conservative. Concurrent, deriving insights into the error for arbitrary observables is non-trivial and still an open problem. Moreover, while the controller performance is optimized regarding the lifted states, the resulting performance of the underlying system might be unclear.

%% file: 07_Feedback_linearization.tex
\section{Data-driven control by approximate nonlinearity cancellation and feedback linearization}\label{Sec_FeedbackLin}

Approximate nonlinearity cancellation aims to find a feedback law that stabilizes the systems while reducing the influence of the nonlinearity of the closed-loop dynamics. Whereas the presented nonlinearity cancellation technique works in the original states, feedback linearization transforms a nonlinear system via a certain change of coordinates and a feedback law \citep{Isidori}. For flat systems, the transformed system dynamics becomes linear as the nonlinear internal dynamics vanishes. Systems of that kind emerges in practice, e.g., in robotics and in automatic flight control \citep{Flat_Systems}. Due to the feedback law in both paradigms, the dynamics is changed in order to obtain a closed-loop description suitable for linear control design techniques. In contrast, the previously presented approaches aim to obtain a suitable characterisation of the open-loop dynamics itself instead of modifying it.  

We report the data-driven control literature exploiting SDPs within these control frameworks. We begin with the data-driven nonlinearity cancellation because this result will be required to obtain a linear system description via feedback linearization.

\subsection{Data-driven approximate nonlinearity cancellation}\label{Sec_NonCan}

In this section, we review exact and approximate nonlinearity cancellation for a data-driven controller design of nonlinear systems by SDPs. More specifically, {\cite{NonlinearityCancell}} aims to stabilize and cancel the nonlinearity of the discrete-time nonlinear system
\begin{equation}\label{Sys_NonCan}
x(t+1)=f(x)+Bu,
\end{equation}
with unknown drift $f:\mathbb{R}^{n_x}\rightarrow\mathbb{R}^{n_x}$ and unknown input matrix $B\in\mathbb{R}^{n_x\times n_u}$. To infer on the nonlinear function $f$ from data, the knowledge of a function basis is supposed.

\begin{assu}[\cite{NonlinearityCancell} (Assumption 1)]\label{Assu_NonCanel}
	Let a continuous function $z:\mathbb{R}^{n_x}\rightarrow\mathbb{R}^{n_z}$ be known such that $f(x)=Az(x)$ for some matrix $A\in\mathbb{R}^{n_x\times n_z}$.
\end{assu}  

A similar assumption is already introduced and clarified in Section~\ref{Sec_LPVEmb} and Section~\ref{Sec_SysRepExLin}. For that reason, a suitable basis $z$ might be available from first principles, while only the system parameters $A$ are unknown. Since one key idea will be to separate the linear and nonlinear dynamics, \cite{NonlinearityCancell} writes $z$ without loss of generality as
\begin{equation*}\label{Sep_linNL}
	z(x)=\begin{bmatrix}x\\q(x)\end{bmatrix},
\end{equation*}
where $q$ only contains nonlinear functions.

After the stage is set, \cite{NonlinearityCancell} extends \cite{PersisLinear} by determining a data-driven representation of the closed loop with system \eqref{Sys_NonCan} and the nonlinear state feedback $u(t)=Kz(t)$. To this end, the data $\{x_i,u_i\}_{i=1}^S$ from system \eqref{Sys_NonCan} is collected into the matrices $U=\begin{bmatrix}u_1 & \cdots & u_{S-1}\end{bmatrix},X^+=\begin{bmatrix}x_2 & \cdots & x_{S}\end{bmatrix},Z=\begin{bmatrix}x_1 & \cdots & x_{S-1}\\ q(x_1) & \cdots & q(x_{S-1})\end{bmatrix}$.
Then, the closed loop can be characterized by the following data-dependent matrices
\begin{align*}
x(t+1)&=\begin{bmatrix}A & B\end{bmatrix} \begin{bmatrix}I\\K\end{bmatrix}z(x)\\
&=	\begin{bmatrix}A & B\end{bmatrix}\begin{bmatrix}Z\\U\end{bmatrix}Gz(x)=X^+Gz(x),
\end{align*}
where the last equality follows from $X^+=AZ+BU$. The second equality holds for a matrix $G$ satisfying 
\begin{equation}\label{Consistency_CL}
\begin{bmatrix}I\\K	\end{bmatrix}=\begin{bmatrix}Z\\U\end{bmatrix}G,
\end{equation} 
as already shown in \cite{PersisLinear} and Section~\ref{Sec_Direct}. Proceeding the idea of separating the linear and nonlinear terms, \cite{NonlinearityCancell} (Lemma 1) derives the data-based closed-loop description $x(t+1)=Lx+Nq(x)$ with $L=X^+G_1$, $N=X^+G_2$, $G=\begin{bmatrix}G_1 & G_2\end{bmatrix}$. 

Given this description of the closed loop, it is globally asymptotically stabilized for a $G$ satisfying $X^+G_2=0$ and that 
renders $L=X^+G_1$ to be Schur. Indeed, by cancelling the input matrix $X^+G_2$ of the nonlinearity, only the linear system matrix $X^+G_1$ needs to become stable. \cite{NonlinearityCancell} (Theorem 1 and Theorem 3) formulates this problem as an SDP. Furthermore, \cite{NonlinearityCancell} (Theorem 4) investigates the case when the nonlinearity $X^+G_2$ can not be cancelled out exactly. Instead, its influence $||X^+G_2||_\text{Fr}$ is minimized. Therefore, the nonlinearity is only approximately cancelled and the origin is only asymptotically stable, if the linear part dominates the nonlinear, i.e.,
\begin{equation*}
\lim\limits_{||x||_2\rightarrow0}\frac{||q(x)||_2}{||x||_2}=0.
\end{equation*}
This procedure can also be adapted for noisy data and bounded neglected nonlinearities \citep{NonlinearityCancell}(Theorem~6 and 8). To this end, the controller again musts render the linear closed-loop dynamics stable. However, since the dynamics can not exactly be identified from data, a set of systems matrices has to be stabilized. 

While \cite{PersisLinear} neglects the nonlinearity and elaborates a linear feedback design, \cite{NonlinearityCancell} achieves a more flexible nonlinear feedback law. Moreover, the effect of the nonlinearity is minimized to increase the region of attraction. However, since the stabilization criterion in \cite{PersisLinear} and \cite{NonlinearityCancell} contains only the linear part of the dynamics, only an asymptotically stable equilibrium can be ensured without guarantees for the size of the region of attraction. Moreover, the refinement of the latter is unclear because the nonlinearity might actually be useful for stabilization. Thus, a cancellation might even decrease the region of attraction or might render the system non-robust regarding small disturbances \citep{Robust_Feedback_lin}. Instead of cancelling the nonlinearity, it is bounded and incorporated into the controller synthesis in \cite{MartinCS} to achieve performance criteria and global stability.

To impose further guarantees, \cite{NonlinearityCancell} proposes to analyze the region of attraction (Proposition 1) or robust invariant sets (Theorem 7) for the obtained closed loop. Nonetheless, here the problem arises that the closed-loop dynamics contains nonlinear terms, and thus estimating these sets calls for solving a nonconvex optimization problem. Furthermore, the complexity of the closed-loop parametrization \eqref{Consistency_CL} increases with the number of samples as in \cite{PersisLinear}. Thus, a set membership for the uncertain system matrices $A$ and $B$ might be preferable. Finally, note that nonlinearity cancellation conceptually does not allow for analyzing the underlying system regarding, e.g., dissipativity.

\subsection{Data-driven feedback linearization}

Before studying the data-based case, we begin with a general recap of feedback linearization of flat single-input single-output discrete-time systems and refer to \cite{Feedback_lin_MIMO} for the multiple-input multiple-output case. Consider the nonlinear system 
\begin{equation}\label{System}
\begin{aligned}
	x(t+1)&=f(x(t),u(t)),\\
	y(t)&=h(x(t)),
\end{aligned}
\end{equation}
with smooth functions $f:\mathbb{R}^{n_x}\times\mathbb{R}\rightarrow\mathbb{R}^{n_x}$ and $h:\mathbb{R}^{n_x}\rightarrow\mathbb{R}$. Further, we assume that system \eqref{System} satisfies
\begin{align*}
	\frac{\partial}{\partial u}(h\circ f_0^i\circ f(x,u))&=0,\quad \forall(x,u)\in\mathbb{R}^{n_x}\times\mathbb{R},0\leq i\leq n_x-2,\\
	\frac{\partial}{\partial u}(h\circ f_0^{n_x-1}\circ f(x,u))&\neq0,\quad \forall(x,u)\in\mathbb{R}^{n_x}\times\mathbb{R},
\end{align*}
with $f_0^i$ denoting the $i$ times composition of $f(x,0)$. Therefore, system \eqref{System} has a well-defined relative degree $n_x$, and hence is called flat. 
Following \cite{Flatness}, one can define the state transformation
\begin{equation*}
	z(t)=\Psi(x(t))=\begin{bmatrix}h(x(t))\\ h\circ f_0(x(t))\\\vdots\\ h\circ f_0^{n_x-1}(x(t))
	\end{bmatrix}=\begin{bmatrix}y(t)\\ y(t+1)\\\vdots\\ y(t+n_x-1)
	\end{bmatrix},
\end{equation*}
which yields the transformed system dynamics of \eqref{System}
\begin{equation}\label{Feddback_System}
\begin{aligned}
	z(t+1)&=\begin{bmatrix}z_2(t)\\ \vdots\\ z_{n_x}(t)\\h\circ f_0^{n_x-1}\circ f(x(t),u(t))
	\end{bmatrix},\\
	y(t)&=z_1(t).
\end{aligned}
\end{equation}
Thus, a flat system can be fully linearized by the state transformation $z=\Psi(x)$, containing only forward time-shifted outputs, and the input transformation $v(t)=h\circ f_0^{n_x-1}\circ f(x(t),u(t))=\Psi_v(x(t),u(t))=\Psi_v(\Psi^{-1}(z(t)),u(t))=\Psi_v(Y(t),u(t))$ with $Y(t)=\begin{bmatrix}y(t)&\cdots&y(t+n_x-1)\end{bmatrix}$.  

After this introduction to feedback linearization, we can proceed to the data-driven control of flat systems. To this end, let the flat system \eqref{System} be unknown, and thus the input transformation $\Psi_v$ as well. Note that the unknown system can be perturbed to check for flatness \citep{FeedbackLin3}. Then {\cite{FeedbackLin2}} suggests to assume that $\Psi_v$ admits a linear combination of known functions.

\begin{assu}[\cite{FeedbackLin2}]\label{Assumption_FeedbackLin}
Suppose that the input transformation takes the form $\Psi_v(Y(t),u(t))=a^T\tilde{\Psi}_v(Y(t),u(t))$ with known vector of functions $\tilde{\Psi}_v:\mathbb{R}^{n_x}\times\mathbb{R}\rightarrow\mathbb{R}^{n_{\tilde{\Psi}}}$ and unknown coefficient vector $a\in\mathbb{R}^{n_{\tilde{\Psi}}}$.
\end{assu}  
In general, finding a suitable choice of $\tilde{\Psi}_v$ can be challenging, in particular for discrete-time systems. Indeed, the function composition $\Psi_v=h\circ f_0^{n_x-1}\circ f(x,u)$ prevents the derivation of $\tilde{\Psi}_v$ from a function basis of the system dynamics $f$, which might be available from first principles. In contrast, this is possible for continuous-time systems as the function composition is replaced by Lie derivatives. For instance, consider the flat discrete-/continuous-time system
\begin{align*}
	x_1(t+1)|\dot{x}_1(t)&=\alpha_1\sin(x_2(t)),\\
	x_2(t+1)|\dot{x}_2(t)&=\alpha_2\cos(x_1(t))+u(t),\\
	y(t)&=x_1(t),
\end{align*}
with unknown coefficients $\alpha_1$ and $\alpha_2$ and basis functions $\sin(x_2),\cos(x_1)$, and $u$. For the discrete-time case, the input transformation corresponds to $\Psi_v(x,u)=\alpha_1\sin(\alpha_2\cos(x_1)+u)$. Hence, $\tilde{\Psi}_v(x,u)=\sin(\alpha_2\cos(x_1)+u)$ which is however unknown due to the unknown coefficient $\alpha_2$. Contrary, in continuous time, $v(t)=\ddot{y}(t)=\alpha_1\alpha_2\cos(x_1)\cos(x_2)+\alpha_1\cos(x_2)u$ with known basis functions $\cos(x_1)\cos(x_2)$ and $\cos(x_2)u$.

In case a function basis $\tilde{\Psi}_v$ is not available and only the scalar product of $\tilde{\Psi}_v$ is required, then \cite{FeedbackLin2} proposes to consider an infinite number of basis functions and to compute the scalar product by kernel methods. However, this comes with the drawback that an infinitely long PE input sequence is required to span the whole input-output trajectories of the system. Alternatively, \cite{FeedbackLin3} (Assumption 5) relaxes Assumption~\ref{Assumption_FeedbackLin} by incorporating an uniformly bounded uncertainty $\epsilon$
\begin{equation}\label{Feedback_Lin_error}
	\Psi_v(Y(t),u(t))=a^T\tilde{\Psi}_v(Y(t),u(t))+\epsilon(Y(t),u(t)).
\end{equation}  
Thereby, $\tilde{\Psi}_v$ does not necessarily have to span a function basis for $\Psi_v$. 

Under Assumption~\ref{Assumption_FeedbackLin}, the feedback linearization \eqref{Feddback_System} of nonlinear system \eqref{System} yields
\begin{equation}\label{Feedback_System}
\begin{aligned}
	z(t+1)&=\underbrace{\begin{bmatrix}0 & 1 &\cdots & 0\\ \vdots & \ddots & \ddots & \vdots\\ \vdots & & \ddots & 1\\ 0 & \cdots & \cdots & 0
	\end{bmatrix}}_{=:D}z(t)+\begin{bmatrix}0\\\vdots\\0\\a^T\end{bmatrix}\tilde{\Psi}_v(Y(t),u(t)),\\
	y(t)&=\begin{bmatrix}1&0&\cdots&0\end{bmatrix}z(t).
\end{aligned}
\end{equation}
Therefore, the behavior $\tilde{\Psi}_v\rightarrow y$ is linear with unknown input matrix $\begin{bmatrix}0&\cdots&0&a\end{bmatrix}^T$. Thus, \cite{FeedbackLin2} (Proposition 1) can apply Willems' fundamental lemma for LTI systems here.
Note that 
the PE condition is more difficult to be fulfilled for a larger function basis. Moreover, the PE condition includes input and output data instead of only input data as in the LTI case. Therefore, the PE condition can only be checked after an experiment is carried out. We refer to \cite{ExpDesign_FeebackLin} (Theorem 5) for the design of inputs to guarantee PE of a sequence of monomial basis functions.

By means of the extension of Willems' fundamental lemma for flat systems, \cite{FeedbackLin2} solves the data-driven simulation and output-matching control problem for feedback linearizable systems. 
\cite{FeedbackLin3} generalises these results for the relaxed version \eqref{Feedback_Lin_error} of Assumption~\ref{Assumption_FeedbackLin}. Specifically, the effect of uniformly bounded uncertainty $\epsilon$ and output measurement noise on the data-driven simulation and output-matching control problem is analyzed by providing output error bounds. Bounding this error together with the extended fundamental lemma can be leveraged for robust data-driven predictive control with rigorous stability guarantees \citep{FeedbackLin,MPC_FeebackLin}.

While the feedback linearization leads to a linear behavior from the new input $v$ to $y$, the system dynamics $u\rightarrow y$ in \eqref{Feedback_System} is nonetheless nonlinear. Therefore, the extended fundamental lemma results in nonlinear optimization problems for the data-driven simulation, output-matching control, and predictive control. Thus, solving these optimization problems might be difficult or computationally expensive, in particular, if a large prediction horizon or number of basis functions are considered. Furthermore, for data $\{x_i,u_i\}_{i=1}^S$ satisfying $x_{i+1}=f(x_i,u_i)$ and error \eqref{Feedback_Lin_error}, we can determine a set membership containing $f(x,u)$
\begin{align*}
	&\left\{\Psi^{-1}(D\Psi(x)+\begin{bmatrix}0\\1\end{bmatrix}(\tilde{a}^T\tilde{\Psi}_v(x,u)+\epsilon(x,u))): ||\epsilon(x,u)||_2\leq \epsilon^*,\right.\\
	&\hspace{0.0cm}\left.\Big|\Big|\Psi(x_{i+1})-D\Psi(x_i)-\begin{bmatrix}0\\1\end{bmatrix}\tilde{a}^T\tilde{\Psi}_v(x_i,u_i)\Big|\Big|_2\leq\epsilon^*,i=1,\dots,S-1\right\}.
\end{align*}
Note that $\Big|\Big|\Psi(x_{i+1})-D\Psi(x_i)-\begin{bmatrix}0\\1\end{bmatrix}\tilde{a}^T\tilde{\Psi}_v(x_i,u_i)\Big|\Big|_2\leq\epsilon^*,i=1,\dots,S-1,$ imply the set of all coefficients $\tilde{a}$ feasible with that data, and thus contains the true unknown coefficients $a$. While the nonlinearity from $\Psi$ and $\Psi^{-1}$ can be circumvented by considering the set membership w.r.t. $z$, the nonlinearity from $\tilde{\Psi}_v(x,u)$ remains.\\

To circumvent a controller design by a nonlinear optimization problem, \cite{NonlinearityCancell} (Assumption 6) restricts Assumption~\ref{Assumption_FeedbackLin} to a linear input transformation concerning the input
\begin{equation}\label{Input_trafo}
	\Psi_v(Y(t),u(t))=a^T\tilde{\Psi}_v(Y(t))+bu(t).
\end{equation}
While this condition might only be possible for special cases, e.g., the polynomial system in
\cite{NonlinearityCancell} (Example 10), it is always satisfied for an input-affine continuous-time system with linear input matrix. For a purely linear system description, \cite{NonlinearityCancell} (Corollary 2) requires in addition to \eqref{Input_trafo} that the nonlinearity $a^T\tilde{\Psi}_v(Y(t))$ can be exactly cancelled. The remaining data-driven inference of the unknowns $a$ and $b$ and the nonlinearity cancellation using SDPs corresponds to the procedure described in Section~\ref{Sec_NonCan}. 

The main advantage of a feedback linearization is that a flat system can exactly be linearized in some coordinates. Thereby, for instance, global asymptotic stability can be ensured. On the other hand, the approximate nonlinearity cancellation from Section~\ref{Sec_NonCan} stays in the original coordinates. Thus, it allows a controller design and closed-loop analysis in the original coordinates. Moreover, the latter does not require flatness of the system.\\

Under the more restrictive assumption that complete dictionaries for the state transformation $\Psi$ and the input transformation $\Psi_v$ are available, \cite{De_Persis_new} learns both transformations from data. To this end, the null space of a data-dependent matrix has to be computed, which is even for large matrices possible. Based on the obtained state and input transformation, an explicit control law to locally stabilize the nonlinear system around an equilibrium can be deduced.\\

In summary, data-driven feedback linearization requires knowledge of a function basis for the input transformation and leads to a nonlinear system description. Hence, it does not allow for a system analysis or state-feedback design by SDPs. A linear description of the system dynamics can only be obtained together with a nonlinearity cancellation. Contrary, Koopman linearization does not require an input transformation. Thus, a (bi-)linear characterization of the system itself is deduced, which is suitable for a system analysis by SDPs. We refer to \cite{Learning_Koopman_Eigenfcn} for a more thorough discussion of advantages of the Koopman paradigm compared to feedback linearization. Moreover, the question rises whether the overhead of SDPs is necessary as practicable and rigorous learning-based approaches based on feedback linearization already exist, e.g., see \cite{Feedbacklin_alternative1} and \cite{Feedbacklin_alternative2}.

%% file: 08_Conclusion.tex
\section{Discussion}\label{Sec_Discussion}

\begin{table*}
	\begin{center}
		\begin{tabular}{c||c|c|c|c}\label{PropertyTable}
			System representation&  Data-driven method type & Application & Prior knowledge  & Noisy data\\\hline\hline
			Polynomial interpolation& Set membership  & SA,SF& First type  & Yes \\
			Polynomial subclass& WFL, DDCLC, set membership & SA,SF& Second type   & Partially  \\
			Data-based closed-loop description& DDCLC& SF& First type   & Yes \\\hline
			Linear sector from GP& Kernel& SA,SF &First type   & Yes\\
			Linearized kernel& Kernel & SA,SF & First type   & Yes\\
			Polynomial kernel& Kernel & SA,SF & First type  &Yes\\\hline
			LPV system& WFL, DDCLC, set membership& SA,PC,SF & Second type   & Partially \\
			LPV embedding& WFL, DDCLC, set membership& SA,PC,SF & Second type   & Partially \\
			Extended linearization&Set membership& PC & Second type   & Yes \\\hline
			Koopman& EDMD& PC,SF & First type   & Partially \\\hline
			Nonlinearity cancellation& DDCLC& SF & Second type   &Yes	\\
			Flat system& WFL& PC & Second type   & Partially	 
		\end{tabular}
		\caption{Properties of the presented data-driven system representations. Abbreviations: system analysis (SA), predictive control (PC), state feedback (SF), Willems' fundamental lemma (WFL),  Data-driven closed-loop characterization (DDCLC).}
	\end{center}	
\end{table*} 

After the presentation of the various data-based system representations, we consider in this section a more general and embraced discussion. Specifically, we discuss the derivation of the system representations and the required prior model knowledge. Moreover, we compare the implementability of the presented approaches with similar data-driven methods from the literature. We review how the different methods incorporate noisy data. Finally, we give some general comments on the advantages of data-based discrete- and continuous-time system representations. Table~2 partially summarizes the discussion.

\subsection{Derivation of data-driven system representations}

\begin{figure}
	\centering
	\begin{tikzpicture}
	\node[inner sep=0pt] (p1) at (0,0) {\includegraphics[width=.5\textwidth]{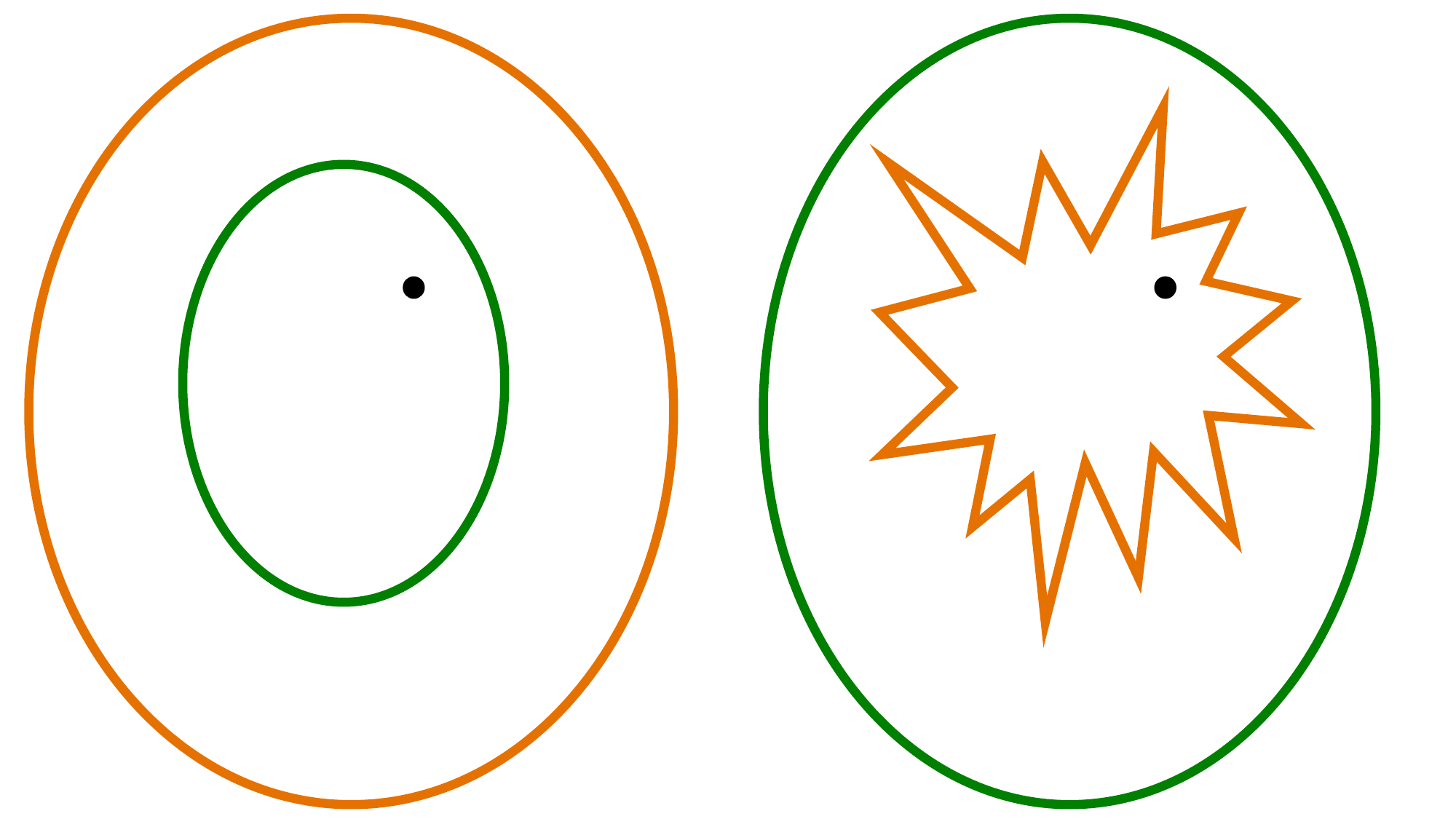}};
	
	\node[inner sep=0pt] (p1) at (-2.5,-0.5) {$\Sigma_\text{\textcolor{green!50!black}{sur}}$};
	\node[inner sep=0pt] (p1) at (-2.8,-1.5) {$\Sigma_\text{\textcolor{green!50!black}{sur}-\textcolor{orange!90!black}{dd}}$};
	\node[inner sep=0pt] (p1) at (-2.6,0.75) {$f(x,u)$};
	
	\node[inner sep=0pt] (p1) at (2,0.0) {$\Sigma_\text{\textcolor{orange!90!black}{dd}}$};
	\node[inner sep=0pt] (p1) at (1.8,-1.8) {$\Sigma_\text{\textcolor{orange!90!black}{dd}-\textcolor{green!50!black}{sur}}$};
	\node[inner sep=0pt] (p1) at (2.13,0.75) {$f(x,u)$};

	\end{tikzpicture}
	\caption{Illustration for the derivation of data-driven system representations for nonlinear systems suitable for system analysis and control by SDPs. Right figure for GP/kernel regression and left figure for the remaining approaches.} 
	\label{Fig.Derivation}
\end{figure}
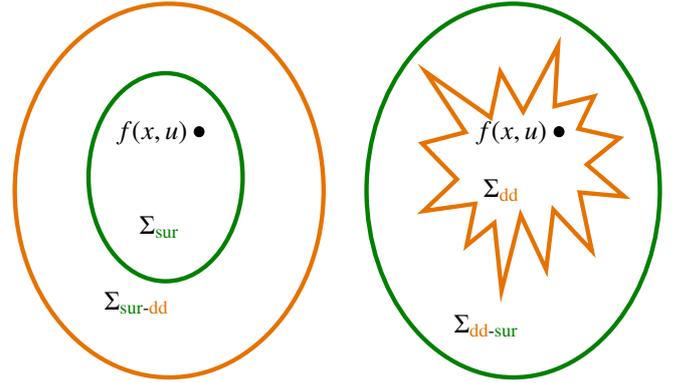

We can identify a common procedure for the derivation of the presented data-driven system representations. Using techniques from model-based control theory, the nonlinear dynamics $f(x,u)$ is first embedded into a set of surrogate systems $\Sigma_\text{sur}$, which are suitable for a system analysis or a controller synthesis via SDP. However, for unknown nonlinear systems, these surrogate systems include unknown coefficients. Therefore, $\Sigma_\text{sur}$ is embedded into the set of systems $\Sigma_\text{sur-dd}$ with data-driven inference on the unknown coefficients. For the polynomial interpolation approach, $\Sigma_\text{sur}$ corresponds to the polynomial sector \eqref{Set_Mem_TP} from TPs and $\Sigma_\text{sur-dd}$ corresponds to \eqref{Set_mem_TP} with the inference on the TP from data. Analogously, for data-driven LPV embedding, $\Sigma_\text{sur}$ is the standard LPV embedding $\{A(p)x+B(p)u:p\in\mathbb{P}\}$ such that $\Sigma_\text{sur-dd}$ corresponds to \eqref{Set_mem_LPV2}. Indeed, we can proceed this for all set membership introduced here except for the ones obtained by GP and kernel ridge regression. There first a data-driven set membership $\Sigma_\text{dd}$ including nonlinear functions is obtained, e.g., \eqref{Set_mem_GP}. Then a set of linear surrogate systems $\Sigma_\text{dd-sur}$, e.g., from \eqref{Com_Sec_bound} is obtained. We illustrate both procedures in Figure~\ref{Fig.Derivation}.

\subsection{Data-driven system analysis, predictive control, and state-feedback design}
The presented system representations can be applied for data-driven system analysis and control. We provide a brief summary in Table~2. Note that the extension of the state feedback with input-state measurements to output feedback with input-output data is possible for a discrete-time setup by considering the extended state of past outputs \citep{PersisLinear, Berberich, LFT_SA}. However, the extended state might lead to conservative requirements on the system and data. \\

\subsection{Prior system knowledge}

All presented system representations require some prior system knowledge to infer guarantees on the dynamics from data. \TM{We divide these into two categories}. 

The \TM{first type of assumption} asks for an upper bound on a Lipschitz constant, higher order partial derivatives, or the complexity of the dynamics. While these \TM{assumptions} are easier to satisfy \TM{than the second type}, they only allow an accurate inference on the dynamics in the neighbourhood of samples. Indeed, \TM{this} information typically implies only a correlation of the nonlinear behaviour at a data point and its neighbourhood. Due to the multitude of local models from polynomial approximation or \TM{the} non-parametric model from kernel regression, these approaches are computational more demanding. This drawback might be circumvented in an online control scheme, where only the local behaviour of the system is relevant at one time step.

The \TM{second type of} assumptions \TM{calls for a function basis that captures the nonlinearity of the dynamics, the scheduling parameter, a state transformation, etc. Thereby, these assumptions} reduce the problem to unknown coefficients. We saw inferences on these coefficients by a set-membership procedure or Willem's fundamental lemma. Due to the prior information on the actual nonlinearity, we will observe in Section~\ref{SecImple} that the number of samples are lower than for the approaches based on the first type of assumption and some machine learning approximation methods, e.g., using neural networks. 

\subsection{Implementability}\label{SecImple}

For the discussion of implementability of the presented approaches, we rely on the examples provided in the reported works. We first observe that numerical examples of a nonlinear system with two states are mostly considered. This is in line with the system complexity of examples of other data-driven control approaches with comparable focus on theoretical guarantees, e.g., backstepping for GPs \citep{GP_Backstepping}, estimation of the region of attraction from GPs \citep{GP_Lyapunov}, control certificate functions based on Lipschitz \TM{estimation} \citep{Lipschitz_online_control}, safe reinforcement learning \citep{Reinf_NL}, neural Lyapunov function \citep{NN_controller}, and learning deep neural networks for Koopman models \citep{NN_control}. These works study the same inverted pendulum example as in \cite{MartinCS}, \cite{LPV_nonlinear}, \cite{Poly_Kernel}, \cite{NonlinearityCancell}, and \cite{Koopman_Guarantees}.

For the inference on a Koopman model, \cite{Koopman_Guarantees} and \cite{Strasser_Koopman} require $5\cdot 10^4$ samples for the inverted pendulum with 5 lifted states and $2000$ for a Van der Pol oscillator with $32$ lifted states, respectively. For learning a neural Lyapunov function or a deep neural network for a Koopman model, \cite{NN_controller} and \cite{NN_control} call for even $10^5$ and $1.6\cdot10^6$ samples, respectively. Contrary, \cite{MartinCS}, \cite{LPV_nonlinear}, \cite{Poly_Kernel}, \cite{NonlinearityCancell}, and \cite{Koopman_Guarantees} require only between $10$ and $100$ samples from the inverted pendulum. As expected, approaches with a priori known nonlinear dictionary call for less data. The computation time to solve the resulting SDPs are reported in \cite{MartinCS} with $8\,\text{s}$, \cite{Poly_Kernel} with $5$ minutes, and \cite{Strasser_Koopman} with less than one second.

This brief comparison shows that the presented data-driven control methods perform in numerical examples comparable with the existing literature, while providing additional advantages as discussed in the introduction. While \cite{MartinTP1}, \cite{LPV_nonlinear}, and \cite{BerberichNL} present initial applications on experimental examples, a broader application of the presented system representations for real data and a more detailed comparison using benchmark examples should be part of future research.\\

\subsection{Noisy data}

As summarized in Table~2, except for the fundamental lemma, most presented system representations provide guarantees though noise-corrupted data. The noise is mostly characterized by deterministic bounds as commonly supposed in data-driven control \citep{vanWaarde}, data-driven system analysis \citep{LFT_SA}, set-membership identification \citep{Milanese}, adaptive control \citep{Adaptive_noise}, and robust model predictive control \citep{Robust_MPC}. However, such a noise description might be conservative if the noise is, e.g., Gaussian distributed as often assumed in system identification. Nevertheless, one can obtain for the presented set memberships, which are linear in the unknown parameters, analogous results with probabilistic guarantees for Gaussian distributed noise following \cite{Umenberger} and \cite{MartinCS}.

Furthermore, an additive measurement noise $d_{\text{meas}}$ on the true states $x_{\text{true}}$, i.e., $x_{\text{meas}}=x_{\text{true}}+d_{\text{meas}}$, yields a more challenging errors-in-variables problem than the presented simplified noise models. To refine these models, one could extend the SOS approach for LTI systems from \cite{Error_in_Var}.\\

\subsection{Discrete- and continuous-time models from data}

Throughout the presentation of the different data-based methods, we have seen continuous-time as well as discrete-time setups. While the approaches from Section~\ref{Sec_PolyApprox}, \ref{Sec_GP}, \ref{Sec_DDLPV} (ii) and (iii), \ref{Sec_Koopman}, and Section~\ref{Sec_NonCan} can be considered for both setups, the question rises what are the advantages compared to each other?

The main advantage of a discrete-time system representation is that only measurements of state or output trajectories are required instead of their time derivatives in addition. The latter can usually only be obtained under weak noise and fast sampling, which allows for signal smoothing with subsequent approximation of the time derivatives by finite differences. 

Advantages of continuous-time system representations are that the learned parameters have a physical meaning as they do not depend on the sampling interval. Moreover, non-uniformly sampled data can be handled. Furthermore, the choice of the sampling interval is less critical than in discrete time. Indeed, if the sampling interval is too low compared to the system dynamics, then the eigenvalues of the discrete-time model are close to the unit circle leading to inaccurate stability inferences. We refer for more details to \cite{CT_SysID}. Lastly, we mention that, e.g., the condition for Lyapunov stability is always linear w.r.t. the dynamics in continuous time while quadratic only for quadratic Lyapunov functions in discrete time. At the same time, transferring these conditions into an SDP by LMI robust control techniques is easier if these conditions are linear or quadratic regarding the nonlinear dynamics.\\

\subsection{A new paradigm of data-driven control methods?}

According to Definition 4 for data-driven control in the survey of \cite{DDC_Survey}, the approaches presented here would not be classified as data-driven. Indeed, the definition requires the controller design to be based on input-output data and does not allow for exploiting any model information. Moreover, instead of a direct method from data to control input, here first a data-based representation of the system is obtained, which is then explicitly included in the controller synthesis by SDPs. Furthermore, in view of \cite{DDC_Survey}, data-driven control methods are applicable independent of the class of systems. Thereby, the paradigm of data-driven control has changed over the last decade to a much broader collection of control techniques.

\section{Conclusion}\label{Sec_Conclusion}

In this survey, we provided an overview on data-driven control approaches for nonlinear systems. In particular, the focus lay on data-based system representations, which are tailored for verifying system properties and designing controllers by SDPs. Thereby, these methods strive to establish a framework for nonlinear data-driven system analysis and control rather than providing specific control schemes.
More specifically, we discussed data-driven control by polynomial approximation, GP and kernel regression, LPV embedding, state lifting, and nonlinearity cancellation and feedback linearization.

Except for GP and kernel regression, these data-driven system representations are inspired by the model-based control literature. There the same techniques are leveraged to derive a verification of system properties and a controller synthesis by SDPs despite nonlinear system dynamics. Since the goal of these control methods is to derive a linear-like representation, data-driven techniques for LTI, bilinear, LPV, and polynomial systems are still relevant in the nonlinear case. Moreover, most of the data-based system characterizations for nonlinear systems that we have presented are combined with LMI-based robust techniques to achieve rigorous guarantees. In contrast, many system representations from system identification and machine learning are tailored to provide a precise surrogate model of the underlying system, and thus aim to approximate its dynamics as precise as possible. However, this typically leads to complex surrogate models preventing a convex system analysis and controller synthesis. 

Guarantees for data-driven inference on nonlinear system regardless of the framework require a priori insights into the dynamics. Hence, the performance of the data-driven approaches not only rely on the informativity of the data but also on the accuracy of the prior information. This is an additional challenge compared to the LTI case, where this kind of assumptions is not required.

To conclude this survey, we provide further open challenges and questions beside the ones mentioned in the individual sections: (i) While the characterization of errors for a polynomial approximation is well-established, the proposed error characterizations, e.g., for the estimation of the Koopman operator, are tailored for control, and thus might be conservative. (ii) How to justify assumptions with prior insights in a data-based setup? (iii) Input-output data are investigated mainly using the extended state \citep{Berberich} resulting in restrictive assumptions on the data. Moreover, \cite{Montenbruck} and \cite{MartinIterative} consider the input-output system behaviour directly by its input-output mapping. However, this requires a large number of measured trajectories. (iv) Closing the gap between representations from a control and a data perspective. While the former is tailored for system analysis and control, the latter for explaining the data and providing a precise surrogate model. However, we search rather for a data-motivated representation suitable for controller design. (v) For a comprising comparison of data-driven control schemes, benchmark systems and data would be required. (vi) The extension of the data-informativity framework \citep{Inform_Survey} to nonlinear systems is missing, i.e., when is the data informative to infer observability, controllability, stabilizability, etc., for a nonlinear system.